\documentclass[11pt]{amsart}
\usepackage[T1]{fontenc}
\usepackage[utf8]{inputenc}
\usepackage{lmodern}
\usepackage{amsmath,amssymb,amsfonts,amsthm,mathrsfs}
\usepackage{graphicx,tikz}
\usepackage{enumitem,booktabs,array}
\usepackage{xcolor}
\usepackage[colorlinks=true]{hyperref}
\usetikzlibrary{decorations.markings}
\hypersetup{pdftitle={Widths of two-dimensional tori and their finite quotients}}

\newtheorem{theorem}{Theorem}[section]
\newtheorem{definition}[theorem]{Definition}

\newtheorem{proposition}[theorem]{Proposition}
\newtheorem{lemma}[theorem]{Lemma}
\newtheorem{remark}[theorem]{Remark}
\newtheorem{corollary}[theorem]{Corollary}
\newtheorem{example}[theorem]{Example}

\newtheorem*{remark*}{Remark}

\numberwithin{equation}{section}

\newcommand{\bM}{\mathbf{M}}

\newcommand{\PP}{\mathbb{P}}

\newcommand{\cA}{\mathcal{A}}

\newcommand{\cE}{\mathcal{E}}

\newcommand{\cG}{\mathcal{G}}
\newcommand{\cH}{\mathcal{H}}
\newcommand{\cI}{\mathcal{I}}

\newcommand{\cK}{\mathcal{K}}

\newcommand{\cR}{\mathcal{R}}

\newcommand{\cZ}{\mathcal{Z}}

\newcommand{\mbfA}{\mathbf{A}}

\newcommand{\mbfF}{\mathbf{F}}

\newcommand{\sH}{\mathscr{H}}

\newcommand{\R}{\mathbb{R}}

\newcommand{\Z}{\mathbb{Z}}
\newcommand{\pr}{\partial}

\newcommand{\scD}{\mathscr{D}}

\newcommand{\scP}{\mathscr{P}}

\newcommand{\scR}{\mathscr{R}}

\newcommand{\Per}{\operatorname{Per}}

\DeclareMathOperator{\area}{Area}

\DeclareMathOperator{\tr}{tr}

\DeclareMathOperator{\sys}{\text{sys}}

\definecolor{quotientglue}{RGB}{43,82,132}
\definecolor{quotientmirror}{RGB}{20,101,91}
\definecolor{quotientcone}{RGB}{164,43,51}
\definecolor{quotientcorner}{RGB}{150,91,0}
\tikzset{
	qcut/.style={draw=quotientglue,line width=.65pt,dash pattern=on 3pt off 1.8pt,
		postaction={decorate},decoration={markings,mark=at position .54 with {\arrow{>}}}},
	qmirror/.style={draw=quotientmirror,line width=1.65pt},
	qcone/.style={circle,fill=quotientcone,draw=quotientcone,inner sep=0pt,minimum size=3.6pt},
	qcorner/.style={rectangle,fill=quotientcorner,draw=quotientcorner,inner sep=0pt,minimum size=4.4pt}
}
\newcommand{\qpt}[4]{\node[#1] at #2 {};\node[#3] at #2 {$#4$};}

\title{Widths of two-dimensional tori and their finite quotients}
\author{Ao Sun}
\date{}

\begin{document}
	
	\begin{abstract}
		We find a combinatorial formula for all the Almgren-Pitts widths of flat tori and their finite isometric quotients. This allows us to compute the widths explicitly for a large class of flat spaces, such as flat tori, squares, equilateral triangles, and right isosceles triangles. One application, which will appear in the future, is the fine asymptotic formula for the volume-spectrum Weyl law for flat tori and their finite isometric quotients up to the fourth term. 
		
	\end{abstract}
	
	\maketitle
	
	\setcounter{tocdepth}{1}
	\tableofcontents
	
	\section{Introduction}
	
	A central idea in geometry for understanding surfaces is to study the collection of curves inside them. For example, the circles of latitude sweep out the round sphere: they start at the south pole, grow to the equator, and shrink to the north pole. The longest curve in this family indicates how ``wide'' the surface is. The minimum among all the lengths of the longest curve in the sweep-out is called the width, showing what the optimal way to measure the wideness of the sphere is. In fact, if we want to pass a rubber band through a sphere, the width is the shortest length that the rubber band needs to have.
	
	We can generalize this idea to a multi-family of curves in a general surface. Roughly speaking, the {\bf $p$-width} $\omega_p$ is the smallest possible length of the longest curve in a topologically nontrivial $p$-parameter family in the space of closed curves in the surface. The sequence $\omega_1\leq\omega_2\leq\cdots$ is called the {\bf volume spectrum}. Gromov suggested viewing it as a nonlinear analogue of the eigenvalues of the Laplacian. In particular, they satisfy a Weyl law similar to the Laplacian eigenvalues.
	
	Due to the nonlinearity and topological complexity, the precise values of widths -- even the first several of them -- are hard to compute. The only known examples are round $2$-spheres computed by Chodosh-Mantoulidis \cite{ChodoshMantoulidis23}, the quotient of the $2$-sphere such as the real projective plane and round hemisphere by Marx-Kuo \cite{MarxKuo25,MarxKuo26}. Even computing the first several widths of a plane polygon is challenging; see \cite{ChodoshCholsaipant25_polygon, Chopra2026computing}.

	The main goal of this paper is to extend the list of surfaces for which we can compute all the widths. Throughout the paper, $\Lambda\subset\R^2$ is a lattice, $T^2=\R^2/\Lambda$ carries the induced flat metric, and $A=\area(T^2)$ is the covolume of $\Lambda$. We write $\omega_p(\Lambda)=\omega_p(T^2)$ for the $p$-width of $T^2$.
	
	Let $\cZ=\cZ_1(T^2;\Z_2)_0$ be the connected component of the space of mod $2$ flat $1$-cycles that contains the zero cycle. By \cite[Sections 2 and 5]{MarquesNeves21}, $\cZ$ consists of the mod $2$ boundaries of sets of finite perimeter. The filling double cover of $\cZ$ associates to each cycle its two complementary fillings. We denote its characteristic class by $\lambda\in H^1(\cZ;\Z_2)$. A $p$-sweepout is a flat-continuous map $\Phi:X\to\cZ$, defined on a finite cubical or simplicial complex $X$, which has no concentration of mass and satisfies $(\Phi^*\lambda)^p\neq0$. The $p$-width is the infimum of $\sup_{x\in X}\bM(\Phi(x))$ over all such $p$-sweepouts. The notion of widths was introduced and studied in many texts, and we only list some of them as references: \cite{Almgren62,Pitts81,Gromov88,Guth09,MarquesNeves17,LMN}. The $p$-widths are closely related to the min-max theory of minimal surfaces. In fact, they are the area of the minimal hypersurfaces constructed by Almgren-Pitts min-max method; see \cite{Pitts81,SchoenSimon81,MarquesNeves16,MarquesNeves17,Zhou20_multi1}. We also highlight that by \cite{GasparGuaraco18,Dey22,ChodoshMantoulidis23}, the widths can also be computed by the singular perturbation method by Allen-Cahn equations; see \cite{HutchinsonTonegawa00,TonegawaWickramasekera12,Guaraco18,GasparGuaraco18,GasparGuaraco19,ChodoshMantoulidis20}.

	Throughout this paper, a convex lattice polygon will always mean a convex polygon with vertices in $\Lambda$, centered at the origin. We also view segments and the point $\{0\}$ as (degenerate) polygons.
	
	\begin{theorem}\label{main_thm_omega_p}
		Let $\Lambda$ be a lattice in $\R^2$ and $T^2=\R^2/\Lambda$. For every integer $p\geq1$, 
		\begin{equation}\label{eq:main}
			\omega_p(T^2)
			=
			\frac12
			\min
			\left\{
			\Per(P):
			\begin{array}{l}
				P\text{ is a centrally symmetric convex lattice polygon},\\
				\#(P\cap\Lambda)\geq p+1
			\end{array}
			\right\}.
		\end{equation}
	\end{theorem}
	To distinguish different tori, we use the notation $\omega_p(\Lambda)$ to denote $\omega_p(\R^2/\Lambda)$.
	
	Note that for any admissible $P_0$ in the right-hand side of \eqref{eq:main}, a potential minimizer has perimeter at most $\Per(P_0)$, so every vertex $v$ satisfies $|v|\leq \Per(P_0)$. Therefore, only finitely many lattice polygons satisfy this bound, and hence the right-hand side of \eqref{eq:main} is the minimum among a finite set. As a result, we can actually compute $\omega_p(\Lambda)$ explicitly by a computer. The lattice-polygon formula also yields an explicit dynamic-programming algorithm that computes all widths up to order $p$ in $O(p^3)$ arithmetic operations.
	
	Because the minimum is taken among symmetric $P$, $\#(P\cap\Lambda)-1$ must be even. Therefore, we can obtain the following corollary immediately.
	\begin{corollary}\label{cor:torus-pairs}
		For every integer $k\geq1$,
		\[
		\omega_{2k-1}(\Lambda)=\omega_{2k}(\Lambda).
		\]
	\end{corollary}

	We next state the formula for general finite isometric quotients of flat tori. Suppose $\Lambda$ is a lattice of $\R^2$ and $T^2=\R^2/\Lambda$ with area $A$. Let $G$ be a finite isometry group of $T^2$ and suppose $Q=T^2/G$. For $g\in G$, we can decompose it by \[
	g(x)=L_gx+t_g\pmod\Lambda
	\]
	where $L_g\in O(2)$ and $t_g$ is a translation. We can define the action of $G$ on $\Lambda$ by
	\[
	g\cdot z=(\det L_g)L_gz.
	\]
	Let us define a special collection of lattice points
	\[
	\Lambda_G
	:=
	\left\{
	z\in\Lambda:
	\frac{\det(z,t_g)}{A}\in\Z
	\text{ whenever }g\cdot z=z
	\right\}.
	\]
	For $z\in\Lambda$, the associated Fourier mode is $e_{Jz/A}(x)=e^{2\pi i\det(z,x)/A}$. Therefore, the condition $z\in\Lambda_G$ means that this mode is invariant under every $g\in G$ satisfying $g\cdot z=z$. In particular, if $z=kz_0$, where $z_0$ is primitive and $k\geq1$, each level set of this mode consists of $k$ parallel closed geodesics in direction $z_0$. The compatibility condition allows these geodesics to be permuted by elements fixing $z$.
	
	For a centrally symmetric, $G$-invariant convex lattice polygon $P$, we define
	\[
	q_G(P)
	=
	\#\bigl((P\cap\Lambda_G)/G\bigr),
	\]
	where the quotient is taken with respect to the action
	$g\cdot z=(\det L_g)L_gz$.  $q_G(P)$ counts the $G$-orbits of compatible frequencies represented by the lattice points of $P$. Equivalently, Lemma \ref{lem:fourier-count}(ii) shows that
	\[
	q_G(P)=\dim_{\mathbb R}V_{JP/A}^G,
	\]
	where $V_{JP/A}^G$ is the space of $G$-invariant real trigonometric polynomials with frequencies in $JP/A$. Note that when counting the dimension of the space, different positive multiples of a primitive lattice vector represent different frequencies, even if they determine the same foliation by parallel geodesics. 
	
	$Q$ may have boundary and interior orbifold points. The widths of $Q$ are defined using sweepouts in the relative mod-$2$ cycle space $\mathcal Z_1(Q,\partial Q;\mathbb Z_2)_0$, with the boundary treated relatively, and the cone points are viewed as interior points.

	\begin{theorem}\label{thm:affine}
		For every integer $p\geq1$,
		\begin{equation}\label{quot:eq:main}
			\omega_p(Q)=\frac1{2|G|}\min
			\left\{
			\Per(P):
			\begin{array}{l}
				P\text{ is a centrally symmetric, $G$-invariant}\\
				\text{convex lattice polygon},\ q_G(P)\geq p+1
			\end{array}
			\right\}.
		\end{equation}
	\end{theorem}
	
	In Appendix \ref{S:quotient classification}, we provide explicit formulas for all seventeen types of compact flat $2$-orbifolds by Conway and Huson \cite{ConwayHuson02}, especially Table II. The spaces include rectangles, the triangles listed in the abstract, and the surfaces of disphenoids. 
	\subsection{Idea of the proof}
	
	The key idea is relating the mass and topology of families of cycles to the perimeter and lattice-point count of lattice convex polygons. The lattice points can be viewed as the homology classes of the torus, which carry the topological information for the cycles; meanwhile the distance between the lattice points can be viewed as the length of the geodesics in the corresponding class, which carry the geometric information of the cycles.
	
	For the upper bound, we use the sweepouts $\partial\{f<0\}$ generated by those $f$ in the linear space of Fourier modes on the torus with frequencies in the dual lattice corresponding to $P$. This construction is suggested by Guth \cite{Guth09} for general spaces, but it works extremely well for the flat tori. We can use Cauchy's perimeter formula, counting the zeros along closed geodesics, to show that
	\[
	\bM(\partial\{f<0\})\leq\frac12\Per(P).
	\]
	This proves the upper bound.
	
	For the lower bound, we associate a centrally symmetric lattice polygon $P_\Sigma$ to each finite-mass cycle $\Sigma$ in the sweepout. It is the convex hull of half the homology classes of all integral representatives of $\Sigma$. Pairing these representatives with constant one-forms and applying Cauchy's perimeter formula again shows that
	\[
	\Per(P_\Sigma)\leq 2\bM(\Sigma).
	\]
	Therefore, the main goal is to show that the topology of a $p$-sweepout forces some associated polygon to contain at least $p+1$ lattice points.
	
	We prove this by constructing equivariant maps on the filling double cover discovered by Marques--Neves in \cite[Section 5]{MarquesNeves21}. For a compact family of cycles with fixed homology data, canonical fillings give a map to $S^0$ when every integral representative is null-homologous. Otherwise, there is a signed relative-area invariant that gives a map to $S^1$. These local maps will be glued by a partition of unity in the image of the sweepouts to get an equivariant map from the space of filling double cover to $S^{p-1}$. As a consequence, this shows that for an $\mbfF$-continuous family $\Phi$:
	\[
	\#(P_{\Phi(x)}\cap\Lambda)\leq p
	\quad\text{for every }x
	\Longrightarrow 
	(\Phi^*\lambda)^p=0.
	\]
	This gives the lower bound via proof by contradiction.
	
	For a finite quotient $Q=T^2/G$, we repeat the above proof but considering $G$-invariant objects lifted from $Q$ to $T^2$. Most of proof is verbatim with $\#(P\cap\Lambda)$ replaced by $q_G(P)$.

	\subsection{Applications}
	The lattice formula has several applications. Here we highlight one that impressed the author the most. As an analogue of the classical Weyl law of the Laplacian operator, the widths also satisfy a Weyl law called the volume spectrum Weyl law. Let $T=\mathbb R^2/\Lambda$ be a flat two-torus, and let $\omega_p(T)$ denote its one-dimensional mod-$2$ volume spectrum. The Weyl law of Gromov--Guth \cite{Gromov88,Guth09} and Liokumovich--Marques--Neves \cite{LMN} shows that there exists a universal constant for all surface $a(1)$ such that
	
	\begin{equation}\label{eq:Weyl}
		\omega_p(T^2)=a(1)\sqrt{\area(T^2) p} +o(\sqrt p).    
	\end{equation}
	Later, Chodosh--Mantoulidis \cite{ChodoshMantoulidis23} showed that the constant $a(1)=\sqrt{\pi}$.
	
	Here we present the higher-order asymptotics, up to order $4$, of the volume spectrum Weyl law. This was discovered by a commercial LLM in the following strange but interesting way: after obtaining the lattice formula, the author asked the LLM what the next-order asymptotics are in the Weyl law, and it could not answer. Then the author asked the LLM to use numerical simulation to test $\omega_p$ of the square tori and Loewner tori, up to $p=5000$. Then the LLM found the next order should be $3.22 p^{-1/6}$. It is very hard to imagine what $3.22$ is, and it turns out to be out of the author's imagination; see \eqref{eq:constants}. Then the author kept raising $p$ and asked the LLM again and again; surprisingly, when $p=10^{18}$, the LLM could only present the numerical estimate but not a definite number; when $p=10^{60}$, the LLM found the precise number with the necessary number-theoretic proof. With this progress, the author asked the LLM to expand the volume spectrum Weyl law to order $4$. The method is quite irrelevant to the current paper, and the author plans to present the proof of these results in the future, in an appropriate way. Here we just present the results.
	
	Let $A=\area(T^2)$ and $\Lambda_0=A^{-1/2}\Lambda$. After a rotation, we may write
	\[
	\Lambda_0=\Lambda_z=y^{-1/2}(\mathbb Z+z\mathbb Z),\qquad z=x+iy\in\mathbb H.
	\]
	The shape $z$ is considered modulo the usual modular action. Let
	\begin{equation}\label{eq:constants}
		\kappa=-\frac{9\,\mathrm B(1/3,1/3)\zeta(-1/3)}{2^{2/3}\pi^{5/6}},\qquad b_z=-\frac{Z_{\Lambda_z}(3/4)}{2\sqrt2\,\pi^{3/4}},
	\end{equation}
	where $\mathrm B$ denotes the beta function and
	\[
	Z_{\Lambda_z}(s)=\sum_{v\in\Lambda_z\setminus\{0\}}|v|^{-2s},\qquad \text{Re}(s)>1,
	\]
	is continued meromorphically elsewhere. Thus $\kappa>0$ is universal, and $b_z$ depends only on the shape.
	
	For a real parameter $n>\pi$, define
	\begin{equation}\label{eq:moving-scales}
		R_n=\sqrt{n/\pi},\qquad \varepsilon_n=R_n^{-8/7},\qquad L_n=R_n^{4/7}\log R_n.
	\end{equation}
	Let $u_j$ be a real orthonormal basis of even Maaß cusp forms on $\operatorname{PSL}_2(\mathbb Z)\backslash\mathbb H$, with respect to $dx\,dy/y^2$, satisfying
	\[
	-y^2(\partial_x^2+\partial_y^2)u_j=(1/4+r_j^2)u_j,\qquad r_j>0.
	\]
	Let
	\[
	E_r(z)=\frac{Z_{\Lambda_z}(1/2+ir)}{2\zeta(1+2ir)}
	\]
	be the Eisenstein series, defined by meromorphic continuation. Write $\rho_j(k)$ and $\rho_r(k)$ for their nonconstant Fourier coefficients, normalized so that
	\[
	u_j(x+iy)=\sqrt y\sum_{k\ne0}\rho_j(k)K_{ir_j}(2\pi|k|y)e^{2\pi ikx},
	\]
	and analogously for the nonconstant part of $E_r$, where $K_\nu$ is the modified Bessel function of the second kind. Let
	\[
	s_r=\frac{1+2ir}{3},\qquad
	I(s)=-\frac{2\zeta(s-1)}{s(1-s)},\qquad
	g_k(s)=2\int_1^\infty[a(a+1)]^{-s}\cos(2\pi ka)\,da
	\quad(k\ge1),
	\]
	where the last integral is understood as an improper integral for $0<\text{Re}(s)<1$. Define
	\[
	\begin{aligned}
		B_j(z)
		&=-\frac{2^{-s_{r_j}}\pi^{1/3-4ir_j/3}\Gamma(ir_j)I(s_{r_j})}{6}\,
		u_j(z)\sum_{k\ne0}g_{|k|}(s_{r_j})|k|^{-ir_j}\overline{\rho_j(k)},\\
		B_{\mathrm{cont}}(r,z)
		&=-\frac{2^{-s_r}\pi^{1/3-4ir/3}\Gamma(ir)I(s_r)}{6}\,
		E_r(z)\sum_{k\ne0}g_{|k|}(s_r)|k|^{-ir}\overline{\rho_r(k)},
	\end{aligned}
	\]
	where $\Gamma$ denotes the gamma function and $B_{\mathrm{cont}}(0,z)$ is defined by continuity.
	
	Finally, define
	\begin{equation}\label{eq:finite-fourth}
		F_n(z)=2\text{Re}\sum_{\substack{j\ \mathrm{even}\\r_j\le L_n}}B_j(z)\exp\left(\frac{ir_j}{3}\log n-\varepsilon_n r_j^2\right)\notag+\frac1{2\pi}\int_{-L_n}^{L_n}B_{\mathrm{cont}}(r,z)\exp\left(\frac{ir}{3}\log n-\varepsilon_n r^2\right)\,dr.
	\end{equation}
	All eigenspaces are included with their full multiplicity. 
	
	\begin{theorem}\label{thm:main_weyl}
		Let $T^2$ be a fixed flat two-torus of area $A$ and shape $z$. Then for every integer $N\ge1$,
		\begin{equation}\label{eq:main_omega_p}
			\omega_p(T^2)=\sqrt{\pi A p}-\kappa\sqrt A\,p^{-1/6}+b_z\sqrt A\,p^{-1/4}+\sqrt A\,p^{-1/3}F_p(z)
			+O_{z,N}\!\left(\sqrt A\,p^{-1/3}\left[p^{-1/42}+(1+\log p)^{-N}\right]\right)
		\end{equation}
		as $p\to\infty$. Moreover,
		\begin{equation}\label{eq:profile-growth}
			F_p(z)=o_z(p^{1/12}).
		\end{equation}
	\end{theorem}
	
	In fact, for every fixed shape $z$, there exist a constant $c_z>0$ and sequences of positive integers $p_k,q_k\to\infty$ such that
	\begin{equation}\label{eq:profile-oscillation}
		F_{p_k}(z)\geq c_z,\qquad F_{q_k}(z)\leq-c_z.
	\end{equation}
	In particular, with the limits taken through the positive integers,
	\[
	\liminf_{p\to\infty}F_p(z)<0<\limsup_{p\to\infty}F_p(z).
	\]

	We next state the quotient formula. Let $Q$ be a compact finite isometric quotient of a flat two-torus. Equivalently, write $Q=\mathbb R^2/\Pi$, where $\Pi$ is a plane crystallographic group. Let $\Lambda$ be its full translation subgroup and let
	\begin{equation}\label{eq:canonical-cover}
		T=\mathbb R^2/\Lambda,\qquad m=[\Pi:\Lambda],\qquad B=\area(Q),\qquad A=\area(T)=mB.
	\end{equation}
	This is the canonical translation cover. Its degree may be smaller than the order of a group in another presentation of the same quotient, because translations must first be absorbed into $\Lambda$. Let $z$ be the shape of $A^{-1/2}\Lambda$, and let $\ell$ be the total length of the mirror boundary of $Q$. Isolated cone points do not contribute to $\ell$. We use the quotient widths in Section \ref{S:widths of quotient}, with relative cycles along the mirror boundary.
	
	\begin{theorem}\label{thm:quotient-transfer}
		For every fixed quotient $Q$ with the notation \eqref{eq:canonical-cover},
		\begin{equation}\label{eq:transfer}
			\omega_p(Q)=\frac1m\omega_{mp}(T)-\frac\ell4+O_Q\!\left(p^{-1/2}\bigl(1+\log\log(p+e^e)\bigr)\right).
		\end{equation}
		The estimate holds at every positive integer index as $p\to\infty$.
	\end{theorem}
	
	\begin{theorem}\label{thm:quotient}
		Under the assumptions of Theorem \ref{thm:quotient-transfer}, for every integer $N\ge1$,
		\begin{equation}\label{eq:quotient-main}
			\omega_p(Q)=\sqrt{\pi Bp}-\frac\ell4-\kappa\sqrt B\,m^{-2/3}p^{-1/6}+b_z\sqrt B\,m^{-3/4}p^{-1/4}
			+\sqrt B\,m^{-5/6}p^{-1/3}F_{mp}(z)+\mathcal R_{Q,N}(p),
		\end{equation}
		where
		\begin{equation}\label{eq:quantitative-error}
			\mathcal R_{Q,N}(p)=O_{Q,N}\!\left(\sqrt B\,p^{-1/3}\left[p^{-1/42}+(1+\log p)^{-N}\right]\right)=o_Q(p^{-1/3}).
		\end{equation}
	\end{theorem}

	\bigskip

	The volume spectrum also determines a flat two-dimensional torus up to isometry; see Theorem \ref{thm:flat_rigidity}. In fact, the first width recovers the systole, the Weyl law recovers the area, and the minimizing lattice polygons determine the remaining parameter.

	\subsection{Organization of the paper}
	We first prove Theorem \ref{main_thm_omega_p}, the width formula for the flat tori. In Section \ref{S:upper bound}, we prove the upper bound of $\omega_p$ in \ref{eq:main}, and in Section \ref{S:lower bound}, we prove the lower bound of $\omega_p$ in \ref{eq:main}, while a major technical part relies on the construction of local equivariant maps will be discussed in Section \ref{S:local equi map}. In Section \ref{S:widths of quotient}, we discuss the extra necessary ingredients in the proof of Theorem \ref{thm:affine}.
	
	We have two appendices: Appendix \ref{S:widths determine a flat torus} discussed the application that the widths determine the flat tori; we leave it in the appendix because the proof is quite tedious. Appendix \ref{S:quotient classification} presents and illustrates the classification of quotient of flat tori by Conway and Huson \cite{ConwayHuson02}.

	\subsection*{Acknowledgement} We are very grateful to Andr\'e Neves. This project and discovery were initiated during the collaboration with him and Fernando Coda Marques \cite{MarquesNevesSun26}, in which paper we realized that the widths of flat tori are related to the Fourier modes of tori. During the preparation of this paper, Andr\'e asked numerous questions and suggested many ideas, which reshaped this project into its current form. We thank Marco Guaraco for his interest, for generously sharing his ideas and some of his current research progress.
	
	\subsection*{AI declaration}
	The Large Language Model (LLM) plays an important role in this project; the author obtained the formulation and proof of the main results via the interaction with a commercial LLM. It also played an essential role in the following mathematical contexts: knowledge of convex geometry, such as zonotopes that are used in Appendix \ref{S:widths determine a flat torus}; the classification of the isometric quotient of flat $2$-tori (recorded in the Appendix \ref{S:quotient classification}); and numerous reference suggestions. We also use an LLM to assist with exposition, proof checking, simplifications, and depicting the figures.

	\section{The upper bound in the main theorem}\label{S:upper bound}
	
	In this section, we prove the upper bound in Theorem \ref{main_thm_omega_p}.
	Following Gromov and Guth \cite{Gromov88,Guth09}, we construct sweepouts from finite-dimensional spaces of trigonometric polynomials. Throughout we fix a lattice $\Lambda\subset\R^2$.
	
	\begin{theorem}\label{thm:torus-upper}
		Let $p\geq1$ be an integer, and let $P=-P$ be a convex lattice polygon such that $\#(P\cap\Lambda)\geq p+1$. Then
		\[
		\omega_p(T^2)\leq\frac12\Per(P).
		\]
	\end{theorem}
	
	\begin{proof}
		Recall that $A=\area(T^2)$. Since
		\[
		\frac{\det(z,\ell)}{A}\in\Z
		\qquad\text{for all }z,\ell\in\Lambda,
		\]
		the functions $x\to e^{2\pi i\det(z,x)/A}$ are well defined on $T^2$ as the Fourier modes. Let
		\[
		V_P=
		\left\{
		f(x)=\sum_{z\in P\cap\Lambda} a_z e^{2\pi i\det(z,x)/A}:  a_{-z}=\overline{a_z}
		\right\}.
		\]
		The Fourier modes with different lattice points are linearly independent. Since $P=-P$, the space $V_P$ consists of real-valued functions and satisfies
		\[
		\dim_{\R}V_P=\#(P\cap\Lambda)\geq p+1.
		\]
		
		For $0\neq f\in V_P$, let $E_f=\{f<0\}$. We first prove that $E_f$ has finite perimeter and
		\[
		\bM(\partial E_f)\leq\frac12\Per(P).
		\]
		Let $J$ denote counterclockwise rotation by $\pi/2$, and write
		\begin{equation}\label{eq:W_K}
			W_P(v)=\max_{z\in P}\langle z,v\rangle
			-\min_{z\in P}\langle z,v\rangle,
			\qquad v\in\R^2.
		\end{equation}
		Fix a primitive vector $w\in\Lambda$, and let $v_w=w/|w|$. The restriction of $f$ to a closed geodesic $t\to x_0+tw$ with $t\in \R/\Z$ is a trigonometric polynomial with integer frequencies $\det(z,w)/A$. Unless it vanishes identically, its number of zeros is at most
		\[
		\frac1A\left(
		\max_{z\in P}\det(z,w)-\min_{z\in P}\det(z,w)
		\right)=\frac{|w|}{A}W_P(Jv_w).
		\]
		To see this, multiply the restriction by a character so that its
		smallest frequency becomes zero, and apply the degree bound for
		polynomials in $e^{2\pi it}$.
		
		Since $f$ is nonzero and real analytic, the slicing characterization of functions of bounded variation \cite{AFP00} shows that $E_f$ has finite perimeter. Let $\nu$ be its measure-theoretic unit normal along $\partial^*E_f$. The transverse circle of the foliation parallel to $w$ has length $A/|w|$, so the slicing formula gives
		\[
		\int_{\partial^*E_f}
		|\langle v_w,\nu\rangle|\,ds
		\leq
		\frac{A}{|w|}\frac{|w|}{A}W_P(Jv_w)
		=W_P(Jv_w).
		\]
		Because the primitive lattice directions are dense in the unit circle, and both sides of the above inequality are continuous in the direction, this inequality holds for every unit vector $v$. Now we integrate this inequality over unoriented unit directions and use Cauchy's perimeter formula \cite{Santalo76,Schneider14},
		
		\begin{equation}\label{eq:fourier-mass}
			2\bM(\partial E_f) \leq
			\int_{S^1/\{\pm1\}}W_P(Jv)\,dv =\Per(P).
		\end{equation}
		Here we used
		$\int_{S^1/\{\pm1\}}|\langle v,\nu\rangle|\,dv=2$.
		
		We also claim that this family has no concentration of mass. In fact, if we choose two linearly independent primitive vectors $w_1,w_2\in\Lambda$, and let $v_i=w_i/|w_i|$. By the same slicing argument as above, for every sufficiently small $\rho$,
		\[
		\int_{\partial^*E_f\cap B_\rho(x)}
		|\langle v_i,\nu\rangle|\,ds
		\leq C_i\rho,
		\qquad i=1,2,
		\]
		where $C_i$ depends only on $P$ and $w_i$. Since $v_1$ and $v_2$ are linearly independent, we have $1\leq C (
		|\langle v_1,\nu\rangle|
		+
		|\langle v_2,\nu\rangle|)$
		for every unit vector $\nu$. Hence
		\[
		\|\partial E_f\|(B_\rho(x))\leq C_P\rho,
		\]
		uniformly in $f$ and $x$. Thus the family has no concentration of mass.
		
		Suppose $V$ is a subspace of $V_P$ with dimension $p+1$. On the unit sphere of $V$, taking the complement sends the continuous map $\Phi(f)=E_f$ to its antipodal involution. Therefore,  the pullback of the filling double cover under $\Phi$ is the antipodal double cover of $\mathbb{RP}^p$. Thus $\Phi^*\lambda$ is the generator of $H^1(\mathbb{RP}^p;\Z_2)$, and $(\Phi^*\lambda)^p\neq0$. Thus, $\Phi$ is a $p$-sweepout. Then by
		\eqref{eq:fourier-mass},
		\[
		\omega_p(T^2)
		\leq
		\sup_{[f]\in\PP V}\bM(\partial E_f)
		\leq
		\frac12\Per(P).
		\]
	\end{proof}
	
	With the theorem, taking the minimum in Theorem \ref{thm:torus-upper} proves the upper bound in \eqref{eq:main}.

	\section{Construction of local equivariant maps}\label{S:local equi map}
	In the following two sections, we will prove the upper bound of Theorem \ref{main_thm_omega_p}. Suppose we are given a $p$-sweepout $\Phi:X\to\cZ$. The lower bound is much more challenging because we need to show there exists a single slice $\Phi(x)$ in the sweepout that has a long enough length, while most of the slices may have short length. We need to use the topology information to show the existence of this long length.

	\subsection{Integral representative}
	
	Let $\Sigma\in\cZ$ have finite mass. An {\bf integral representative of $\Sigma$ modulo $2$} (or {\bf integral representative} for short) is an integral $1$-current $U$ in $T^2$ such that
	\[
	\partial U=0,
	\qquad
	U\bmod2=\Sigma,
	\qquad
	\bM(U)=\bM(\Sigma).
	\]
	We denote the set of such representatives by $\cI(\Sigma)$. For example, if $\Sigma$ is a finite disjoint union of simple closed Lipschitz curves, an integral representative of $\Sigma$ is a choice of orientations on each curve.
	
	The set $\cI(\Sigma)$ is nonempty. Indeed, by \cite{MarquesNeves17}, $\Sigma=\pr U$ for a Caccioppoli set $U$, and the boundary of the integral current $[[U]]$ is an integral representative of $\Sigma$. Since $\Sigma\in\cZ$ is null-homologous with $\Z_2$ coefficients, if $\tilde\Sigma\in \cI(\Sigma)$, then the homology class $[\tilde\Sigma]\in H_1(T^2;\Z)\cong\Lambda$ of every integral representative lies in $2\Lambda$.
	
	We first collect the homology classes represented by integral representatives. For $\Sigma\in\cZ$ and $z\in\Lambda$, we define
	\[
	\cR(\Sigma)=\{[U]/2:U\in\cI(\Sigma)\}\backslash\{0\},
	\qquad
	\cI_z(\Sigma)=\{U\in\cI(\Sigma):[U]=2z\}.
	\]
	The set $\cR(\Sigma)$ is finite and invariant under $x\to -x$. 
	
	\begin{definition}\label{def:asso polygon}
		The {\bf associated polygon} $P_\Sigma$ of $\Sigma$ is the convex hull in $\R^2$ of the set $\{[\tilde\Sigma]/2:\tilde\Sigma\in \cI(\Sigma)\}$.
	\end{definition}
	
	Since the boundaries of both fillings are null-homologous integral representatives, $P_\Sigma$ is the convex hull of $\cR(\Sigma)\cup\{0\}$. In particular,
	\begin{equation}\label{eq:latticecount}
		1+\#\cR(\Sigma)\leq\#(P_\Sigma\cap\Lambda).
	\end{equation}
	
	Pairing with constant $1$-forms shows that $|[\tilde\Sigma]|\leq\bM(\Sigma)$ for every $\tilde\Sigma\in\cI(\Sigma)$. Thus the integral representatives of $\Sigma$ realize only finitely many homology classes, so $P_\Sigma$ is a convex lattice polygon. Since $-\tilde\Sigma\in\cI(\Sigma)$ whenever $\tilde\Sigma\in\cI(\Sigma)$, the polygon $P_\Sigma$ is centrally symmetric with respect to the origin.
	
	Now we state the main result that we will prove in this section. The following proposition gives two local equivariant maps.
	
	\begin{proposition}\label{prop_equiv_maps}\label{prop:local equiv map}
		Let $C\subset\cZ$ be an $\mbfF$-compact family of bounded mass, and suppose that $\cR(\Sigma)=S$ for all $\Sigma\in C$.
		\begin{enumerate}[label=(\roman*)]
			\item If $S=\emptyset$, then the cycles in $C$ admit a continuous choice of filling. Equivalently, the filling double cover over $C$ admits a continuous equivariant map to $S^0\cong\{\pm1\}$.
			\item If $S\neq\emptyset$, then the filling double cover over $C$ admits a continuous equivariant map to $S^1$.
		\end{enumerate}
	\end{proposition}
	
	We prove this proposition in the rest of this section. There is a figure illustrating the proof at the end of this paper; see Figure \ref{fig:local-equivariant-maps}.

	\subsection{Inessential sets and their properties}\label{SS:ineseential prop}
	
	\begin{definition}
		We say that $\Sigma$ is {\bf homologically inessential}, or {\bf inessential} for short, if $\cR(\Sigma)=\emptyset$, that is, if every integral representative of $\Sigma$ is null-homologous.
	\end{definition}
	
	The inessential cycles are important, not only because they are the cycle appearing in Proposition \ref{prop:local equiv map} item (i), but also because they contribute to the difference of two integral representatives of a general cycle. We also define the Borel set to be inessential.

	\begin{definition}\label{def:inessential}
		For a Borel set $\sH\subset T^2$ with $\cH^1(\sH)<\infty$, let $\cG(\sH)$ be the group of integral $1$-cycles $Z$ such that $\|Z\|\leq N\cH^1\llcorner\sH$ for some $N<\infty$. We say that $\sH$ is {\bf inessential} if every cycle in $\cG(\sH)$ is null-homologous.
	\end{definition}
	
	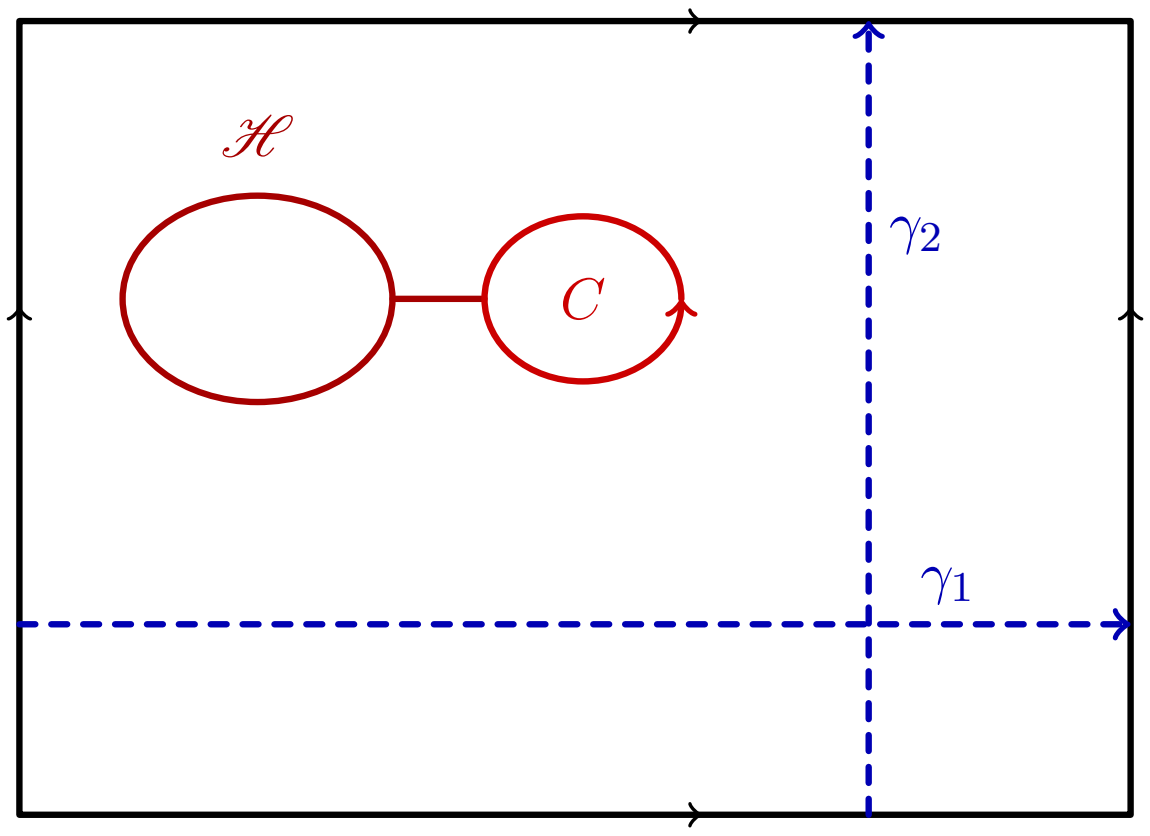
\begin{figure}[ht]
		\centering
		\begin{tikzpicture}[scale=1.05,line cap=round,line join=round]
			
			\draw[very thick] (0,0) rectangle (7,5);
			
			\draw[->,thick] (2.7,0) -- (4.3,0);
			\draw[->,thick] (2.7,5) -- (4.3,5);
			
			\draw[->,thick] (0,1.8) -- (0,3.2);
			\draw[->,thick] (7,1.8) -- (7,3.2);
			
			\node[below] at (3.5,-0.35)
			{$T^2=\mathbb R^2/\Lambda$};
			
			\draw[red!65!black,very thick]
			(1.50,3.25) ellipse [x radius=0.85,y radius=0.65];
			
			\draw[red!65!black,very thick]
			(2.35,3.25) -- (2.93,3.25);
			
			\node[red!65!black,above] at (1.50,4.02) {$\sH$};
			
			\draw[red!80!black,very thick,->]
			(4.17,3.25)
			arc[start angle=0,end angle=360,
			x radius=0.62,y radius=0.52];
			
			\node[red!80!black] at (3.55,3.25) {$C$};
			
			\draw[blue!70!black,very thick,dashed,->]
			(0,1.20) -- (7,1.20);
			
			\node[blue!70!black,above] at (5.85,1.20)
			{$\gamma_1$};
			
			\draw[blue!70!black,very thick,dashed,->]
			(5.35,0) -- (5.35,5);
			
			\node[blue!70!black,right] at (5.35,3.65)
			{$\gamma_2$};
			
		\end{tikzpicture}
		
		\caption{An inessential set $\sH\subset T^2$, drawn as a
			finite graph with smooth edges contained in an embedded disk. In particular, $\mathcal H^1(\sH)<\infty$, and every integral $1$-cycle in $\cG(\sH)$, such as $C$, is null-homologous in $T^2$. The dashed closed geodesics $\gamma_1$ and $\gamma_2$ represent a basis of $H_1(T^2;\Z)$; neither belongs to $\cG(\sH)$.}
		\label{fig:inessential-set}
	\end{figure}
	
	For each finite mass cycle $\Sigma\in\cZ$, we can identify it with a countably $1$-rectifiable Borel set $\|\Sigma\|$ satisfying $\|\Sigma\|=\cH^1\llcorner\Sigma$. Such set representatives are understood up to $\cH^1$-zero measure sets. Let $\pi:\R^2\to T^2$ be the covering map. A {\bf simple Lipschitz loop} is an injective Lipschitz map from a circle, regarded as an integral cycle of multiplicity one. We orient $T^2$ by the standard orientation of $\R^2$ and write $dA$ for its area form. Every integral $2$-current on $T^2$ has the form $\Omega=\theta_\Omega[[ T^2]]$, where $\theta_\Omega$ is an integer-valued function of bounded variation, and its signed area satisfies
	\begin{equation}\label{eq:pairing}
		\Omega(dA)=\int_{T^2}\theta_\Omega\,dA,
		\qquad |\Omega(dA)|\leq\bM(\Omega).
	\end{equation}
	
	The next theorem bounds the homology classes of integral cycles whose mass measures are dominated by a multiple of $\|\Sigma\|$.
	
	\begin{theorem}\label{thm:relax}
		Let $\Sigma\in\cZ$ with $\|\Sigma\|\neq 0$. For every integer $n\ge1$ and every integral $1$-cycle $T$ satisfying $\|T\|\leq n\|\Sigma\|$ as measures, we have $\frac{[T]}{2n}\in P_\Sigma$.
	\end{theorem}
	
	\begin{proof}
		Let us fix a linear functional $\ell$ on $\R^2$, and choose $U\in\cI(\Sigma)$ maximizing $\ell([U]/2)$. Such an integral representative exists because we only need to examine finitely many homology classes.
		
		Consider the integral $1$-cycle $R=nU-T$. Since $\|T\|\leq n\|\Sigma\|$, the multiplicity of $T$ relative to the integral representative of $U$ is an integer in $[-n,n]$ almost everywhere. Hence $R$ has nonnegative multiplicity relative to this integral representative. By the decomposition theorem for integral $1$-cycles \cite{Federer69,BDP}, we can write
		\[
		R=\sum_j C_j,\qquad \sum_j\bM(C_j)=\bM(R),
		\]
		where the $C_j$ are simple Lipschitz loops, and the series converges in mass. Mass additivity excludes cancellation, so each $C_j$ satisfies $|C_j|\leq|\Sigma|$ and agrees with the integral representative of $U$ almost everywhere on $C_j$.
		
		It follows that $U-2C_j$ is again an integral representative of $\Sigma$: it reverses $U$ along $C_j$ and coincides with $U$ elsewhere. By the choice of $U$,
		\[
		\ell\bigl([U]/2-[C_j]\bigr)\le\ell([U]/2),
		\]
		so $\ell([C_j])\geq 0$. Since the decomposition converges in mass, we may evaluate it on the constant closed $1$-form corresponding to $\ell$, which gives
		\[
		n\ell([U])-\ell([T])
		=\ell([R])=\sum_j\ell([C_j])\ge0.
		\]
		Dividing by $2n$, we obtain $\ell([T]/(2n))\le\ell([U]/2)=\max_{P_\Sigma}\ell$. Since $\ell$ is arbitrary and $P_\Sigma$ is convex, we conclude that $\frac{[T]}{2n}\in P_\Sigma$.
	\end{proof}
	
	The following lemma characterizes two classes of inessential sets. Given $\Sigma\in\cZ$, we fix a measurable unit tangent field $\tau$ along $\Sigma$. Every integral representative $U$ of $\Sigma$ can then be written as $U=\sigma_U\tau\|\Sigma\|$, where $\sigma_U:\Sigma\to\{\pm1\}$ is defined almost everywhere.
	
	\begin{lemma}\label{lem:iness}
		Let $\Sigma\in\cZ$.
		\begin{enumerate}
			\item[(i)] If $\Sigma$ is inessential, then every integral $1$-cycle $Z$ satisfying $\|Z\|\leq N\|\Sigma\|$ for some $N<\infty$ is null-homologous.
			\item[(ii)] If $z$ is an extreme point of $P_\Sigma$ and $U_1,\cdots,U_m\in\cI_z(\Sigma)$, then the orientation-disagreement set
			\[
			\scD(U_1,\cdots,U_m):=\{x\in\Sigma:\sigma_{U_1}(x),\cdots,\sigma_{U_m}(x)
			\text{ are not all equal}\}
			\]
			is inessential.
		\end{enumerate}
	\end{lemma}
	
	Roughly speaking, item (ii) implies that for an extreme class, all orientation disagreement is concentrated on an inessential set.
	
	\begin{proof}
		By the mass-additive decomposition into simple Lipschitz loops, it suffices to show that every simple Lipschitz loop whose mass measure is concentrated on the set in question is null-homologous.
		
		(i) Let $C$ be a simple Lipschitz loop whose mass measure is concentrated on $\Sigma$. Since $C$ has multiplicity one, $\|C\|\leq\|\Sigma\|$. Theorem \ref{thm:relax} with $n=1$ yields $[C]/2\in P_\Sigma=\{0\}$, and thus $[C]=0$.
		
		(ii) If $m=1$, the disagreement set is empty, and the statement is automatically true. Suppose that $m\geq 2$, and let $C$ be a simple Lipschitz loop whose mass measure is concentrated on $\mathscr D(U_1,\cdots,U_m)$. Consider the integral cycles
		\[
		T_+=\sum_{i=1}^m U_i+C,
		\qquad
		T_-=\sum_{i=1}^m U_i-C.
		\]
		On the disagreement set, at least one of the signs $\sigma_{U_i}$ differs from the others, so
		\[
		\left|\sum_{i=1}^m\sigma_{U_i}\right|\le m-2.
		\]
		Therefore, either adding or subtracting $C$ yields absolute multiplicity at most $m-1$ there. Outside the disagreement set, $C$ vanishes and the absolute multiplicity is at most $m$. Hence $\|T_\pm\|\leq m\|\Sigma\|$, and Theorem \ref{thm:relax} yields
		\[
		z+\frac{[C]}{2m}\in P_\Sigma,
		\qquad
		z-\frac{[C]}{2m}\in P_\Sigma.
		\]
		The midpoint of these two points is $z$. Since $z$ is an extreme point of $P_\Sigma$, both points coincide with $z$, and therefore $[C]=0$.
	\end{proof}
	
	\subsection{Case (i) of Proposition \ref{prop_equiv_maps}}\label{SS:inessential case}
	
	In order to construct a continuous choice of filling, we need to locally provide a canonical choice of the filling. We do it by picking the filling associated to the closed region bounded by a Jordan curve in the universal cover.
	
	Throughout this subsection, $\sH$ is an inessential set. Every simple Lipschitz loop $C\in\cG(\sH)$ is null-homologous and therefore null-homotopic. Therefore, each lift of $C$ to the universal cover is a planar Jordan curve. Let $[[D_C]]$ be the integral current given by the bounded component of the complement of this lift, oriented so that its boundary is the lift of $C$. Its pushforward to $T^2$ does not depend on the choice of the lift, and the isoperimetric inequality \cite{Osserman78} gives
	\begin{equation}\label{eq:isop-loop}
		\partial\pi_\#[[ D_C]]=C,\qquad
		\bM(\pi_\#[[ D_C]])
		\leq\frac{\bigl(\cH^1(C)\bigr)^2}{4\pi}.
	\end{equation}
	We will choose these fillings compatibly with addition of cycles.
	
	\begin{lemma}\label{lem:coherent}
		Let $(C_i)_{i\geq 1}$ be an at most countable family of simple Lipschitz loops in $\cG(\sH)$ with $\sum_i\cH^1(C_i)<\infty$. Then the lifts of the loops $C_i$ to $\R^2$ can be chosen so that, for $\cH^1$-almost every point $x$ of their union, all loops passing through $x$ use the same lift of $x$. Consequently, for any signs $\varepsilon_i=\pm1$,
		\begin{equation}\label{eq:zerorel}
			\sum_i\varepsilon_i C_i=0
			\quad\to\quad
			\sum_i\varepsilon_i\pi_\#[[ D_{C_i}]]=0.
		\end{equation}
	\end{lemma}
	
	\begin{proof}
		Form a graph with one vertex for each loop $C_i$, joining two vertices when the corresponding loops intersect. In each connected component of this graph, we fix the lift of one loop and propagate the choice of lifts along the edges by requiring that the lifts of adjacent loops agree at an intersection point. A closed path in this graph may change the chosen lift by a deck translation. The arcs of the loops traversed along such a closed path form an integral cycle in $\cG(\sH)$ whose homology class is this deck translation. Since $\sH$ is inessential, this translation vanishes and the lifts are well defined. Applying the same argument to two intersection points of a pair of loops shows that the lifts agree at every common point.
		
		Suppose that $\sum_i\varepsilon_iC_i=0$. Coherence of the lifts implies that the lifted and projected sums have the same multiplicities, so the lifted sum also vanishes. By \eqref{eq:isop-loop}, the series $\sum_i\varepsilon_i[[ D_{C_i}]]$ converges in mass, and its boundary is the lifted sum, which is zero. By the constancy theorem \cite{Federer69,Simon83}, a $2$-current in $\R^2$ with finite mass and zero boundary vanishes. Pushing forward to $T^2$, we obtain \eqref{eq:zerorel}.
	\end{proof}
	
	We now construct a canonical filling on the whole group $\cG(\sH)$. Note that $\cG(\sH)$ is more than what we need for (i) of Proposition \ref{prop:local equiv map}, but this generality is needed to fill in inessential parts of essential cycles for the proof of case (ii).
	
	\begin{proposition}\label{prop:A}
		Let $\sH\subset T^2$ be inessential. There is an additive map $Z\to\cA(Z)$ from the finite-mass integral $1$-cycles whose mass measures are concentrated on $\sH$ to the integral $2$-currents on $T^2$ such that
		\[
		\partial\cA(Z)=Z,
		\qquad
		\bM(\cA(Z))\leq\frac{\bM(Z)^2}{4\pi}.
		\]
		Moreover, for every simple Lipschitz loop $C$ supported on $\sH$, $\cA(C)=\pi_\#[[D_C]]$. In addition, $\cA(Z)$ depends only on $Z$ and not on the choice of the inessential set $\sH$ containing its support.
		
		With this $\cA$, we define
		\begin{equation}\label{eq:alpha}
			\alpha(Z)=\cA(Z)(dA)/A.
		\end{equation}
		Then
		\[
		|\alpha(Z)|\leq\bM(Z)^2/(4\pi A),
		\]
		and $\alpha$ is additive on cycles whose mass measures are concentrated on a common inessential set.
	\end{proposition}
	
	The number $\alpha(Z)$ is the normalized signed area of the canonical filling $\cA(Z)$ of $Z$. Any other integral filling of $Z$ differs from $\cA(Z)$ by an integral multiple of $[[T^2]]$, so its normalized signed area differs from $\alpha(Z)$ by an integer. The canonical choice $\cA$ removes this ambiguity, which is essential in Section \ref{SS:case2}.
	
	\begin{proof}
		By the decomposition theorem for integral $1$-cycles \cite{Federer69,BDP}, we can write
		\[
		Z=\sum_i C_i,
		\qquad
		\sum_i\cH^1(C_i)=\bM(Z),
		\]
		where the $C_i$ are simple Lipschitz loops. Since $\|C_i\|\leq\|Z\|$, each $C_i$ is supported on $\sH$. Define
		\[
		\cA(Z):=\sum_i\pi_\#[[D_{C_i}]].
		\]
		By \eqref{eq:isop-loop}, this series converges in mass and
		\[
		\partial\cA(Z)=Z,
		\qquad
		\bM(\cA(Z))
		\leq\frac1{4\pi}\sum_i\cH^1(C_i)^2
		\leq\frac{\bM(Z)^2}{4\pi}.
		\]
		
		By Lemma \ref{lem:coherent}, $\cA(Z)$ is independent of the decomposition. Lemma \ref{lem:coherent} can also be applied to decompositions of $Z$, $Z'$, and $Z+Z'$ to show the additivity. 
		
		Finally, if $\cA'$ is another map with the stated properties and $R_N=Z-\sum_{i\leq N}C_i$, then $\bM(R_N)\to0$, and hence $\bM (\cA'(R_N) )\leq \bM(R_N)^2/(4\pi)\to0$.
		Since $\cA'(C_i)=\pi_\#[[D_{C_i}]]$, this implies $\cA'=\cA$. The assertions concerning $\alpha$ follow from \eqref{eq:pairing}.
	\end{proof}
	
	The following proposition shows that the map $\cA$ is continuous. In particular, we won't get an extra copy of $[[T^2]]$ in the limit.
	
	\begin{proposition}\label{prop:Acont}
		Let $\sH_i$ and $\sH$ be inessential sets, and let $Z_i\in\cG(\sH_i)$ and $Z\in\cG(\sH)$ be integral cycles of multiplicity one. If $Z_i\to Z$ in the flat topology and $\bM(Z_i)\to\bM(Z)$, then $\cA(Z_i)\to\cA(Z)$ in $L^1({T^2})$.
	\end{proposition}
	
	\begin{proof}
		Let us choose a decomposition of simple Lipschitz loops $Z_i=\sum_j C_{i,j}$, and lift each $C_{i,j}$ to a loop $\widetilde C_{i,j}\subset\R^2$ intersecting a fixed fundamental domain. Let $S_i$ be the sum of their signed Jordan fillings. Then
		\[
		\partial S_i=\sum_j\widetilde C_{i,j},\qquad
		\pi_\#S_i=\cA(Z_i),\qquad
		\bM(\partial S_i)=\bM(Z_i).
		\]
		Since $\bM(Z_i)$ is uniformly bounded, the supports of $S_i$ lie in a fixed compact subset of $\R^2$, and the isoperimetric inequality gives a uniform bound for $\bM(S_i)$. By compactness, after passing to a subsequence, we have $S_i\to S$ in $L^1$.
		Moreover,
		\[
		\pi_\#\partial S=Z,\qquad
		\bM(Z)\leq\bM(\partial S)
		\leq\liminf_i\bM(Z_i)=\bM(Z).
		\]
		Hence there is no cancellation under projection.
		
		Now we decompose $\partial S$ into planar simple Lipschitz loops, where each such loop projects injectively to $T^2$. In fact, if there is an arc joining two points different by a translation in the universal cover, it would project to a cycle with nonzero homology whose mass measure is dominated by $\|Z\|$, which contradicts the inessentiality of $\sH$. Thus their projections form a mass-additive decomposition of $Z$. By the constancy theorem and the definition of $\cA$, we have $\pi_\#S=\cA(Z)$.
		Therefore, applying the argument above to all fillings $S_i$, we have
		\[
		\cA(Z_i)=\pi_\#S_i\to\pi_\#S=\cA(Z)
		\]
		in $L^1$. Since the same argument applies to every subsequence, the whole sequence converges.
	\end{proof}
	
	Let us present a lemma in geometric measure theory. The following lemma follows from the Federer--Fleming compactness theorem \cite{FedererFleming60}.
	
	\begin{lemma}\label{lem:closedgraph}
		Suppose that $\Sigma_i\to \Sigma$ in the $\mbfF$-topology and that $U_i\in \cI(\Sigma_i)$ converge in the flat topology to $U$. Then $U\in \cI(\Sigma)$, and $[U_i]=[U]$ for all sufficiently large $i$. Consequently, over every $\mbfF$-compact family of bounded mass, the space of integral representatives and each subspace with fixed homology class are compact.
	\end{lemma}
	
	\begin{proof}
		Flat convergence gives $\partial U=0$ and $U\bmod2=\Sigma$. Since $\mbfF$-convergence implies convergence of masses, we have
		\[
		\bM(\Sigma)\leq\bM(U)\leq\liminf_{i\to\infty}\bM(U_i)=\lim_{i\to\infty}\bM(\Sigma_i)=\bM(\Sigma),
		\]
		so $U\in \cI(\Sigma)$. Pairing with constant closed $1$-forms shows that $[U_i]\to[U]$. Since these classes lie in the discrete set $\Lambda$, we obtain $[U_i]=[U]$ for all sufficiently large $i$. The compactness assertion follows from the Federer--Fleming compactness theorem \cite{FedererFleming60}.
	\end{proof}
	
	\begin{proof}[Proof of Proposition \ref{prop_equiv_maps}(i)]
		For an inessential $\Sigma$ and $U\in \cI(\Sigma)$, let us define $B(\Sigma)=\cA(U)\bmod2$. By Lemma \ref{lem:iness} (i), $\cA(U)$ is well defined, and $B(\Sigma)$ is a filling of $\Sigma$. We show that $B(\Sigma)$ does not depend on the choice of $U$. In fact, if $U'\in \cI(\Sigma)$, then the integral cycle $(U-U')/2$ satisfies $\|(U-U')/2\|\leq\|\Sigma\|$, and hence belongs to $\cG(\Sigma)$. Then the additivity implies that
		\[
		\cA(U)-\cA(U')=2\cA((U-U')/2).
		\]
		Thus $B(\Sigma)$ does not depend on the choice of $U$.
		
		Let $\Sigma_i\to \Sigma$ under the $\mbfF$-topology in $C$, and choose $U_i\in \cI(\Sigma_i)$. By Lemma \ref{lem:closedgraph}, every subsequence has a further subsequence converging in the flat topology to some $U\in \cI(\Sigma)$. Since $\bM(U_i)=\bM(\Sigma_i)\to\bM(\Sigma)=\bM(U)$, Proposition \ref{prop:Acont} shows that $B(\Sigma_i)\to B(\Sigma)$ in $L^1$. Hence $B$ is continuous. Assigning $1$ to $B(\Sigma)$ and $-1$ to its complement gives the desired map to $S^0\cong\{\pm 1\}$.
	\end{proof}
	
	\subsection{Case (ii) of Proposition \ref{prop_equiv_maps}}\label{SS:case2}
	
	For case (ii) of Proposition \ref{prop_equiv_maps}, intuitively, we define a relative area (modulo the area of $T^2$) to each filling. Then, the relative area gives a map to $S^1$, but it is not enough due to the subtle scenario; see Example \ref{ex:Theta}. Therefore, we need to modify it.
	
	Fix $z\in\Lambda\backslash\{0\}$ and a fixed reference integral cycle $R_z$ of homology class $z$. For $U\in \cI_z(\Sigma)$, the cycle $U-2R_z$ is null-homologous. If $\partial\Omega=U-2R_z$, define
	\begin{equation}\label{eq:Xi}
		\mbfA(U)=\exp\bigl(2\pi i\,\Omega(dA)/A\bigr).
	\end{equation}
	Any two choices of $\Omega$ differ by an integral multiple of $[[{T^2}]]$, so $\mbfA(U)$ is well defined. To obtain an equivariant map, we refine $\mbfA$ as follows.
	
	\begin{definition}\label{def:Theta}
		For $U\in \cI_z(\Sigma)$ and a filling $E$ of $\Sigma$, we choose an integral $2$-current $\Omega$ such that
		\begin{equation}\label{eq:Omega}
			\partial\Omega=U-2R_z,\qquad \Omega\bmod2=E,
		\end{equation}
		and define the {\bf area invariant} $\Theta(U,E)=\exp(\pi i\,\Omega(dA)/A)$.
	\end{definition}
	
	Such an $\Omega$ exists. In fact, the null-homologous cycle $U-2R_z$ bounds an integral $2$-current, and after mod $2$ it becomes a filling of $\Sigma$, hence it is either $E$ or $E^c$. Any two currents satisfying \eqref{eq:Omega} differ by an even multiple of $[[{T^2}]]$, so $\Theta$ is well defined. Moreover,
	\begin{equation}\label{eq:odd}
		\Theta(U,E^c)=-\Theta(U,E),
	\end{equation}
	and $\Theta(U,E)^2=\mbfA(U)$, where $\mbfA$ is defined in \eqref{eq:Xi}.
	
	For an extreme point $z$ of $P_\Sigma$ and $U,V\in \cI_z(\Sigma)$, we define
	\begin{equation}\label{eq:ZUV}
		Z_{UV}=\frac12(U-V).
	\end{equation}
	We have $Z_{UV}\in\cG(\scD(U,V))$, where $\scD(U,V)$ is the orientation-disagreement set. Since $\scD(U,V)$ is inessential by Lemma \ref{lem:iness} (ii), the number $\alpha(Z_{UV})$ is well-defined as in \eqref{eq:alpha}. Moreover, there is a cocycle condition
	\begin{equation}\label{eq:ZUVadd}
		Z_{UV}+Z_{VW}=Z_{UW}.
	\end{equation}
	
	\begin{lemma}\label{lem:Theta}
		The area invariant has the following properties.
		\begin{enumerate}
			\item[(i)] $\Theta$ is continuous with respect to $\mbfF$-convergence of $\Sigma$, flat convergence of $U$, and $L^1$ convergence of $E$.
			\item[(ii)] If $z$ is an extreme point of $P_\Sigma$ and $U,V\in \cI_z(\Sigma)$, then
			\[
			\Theta(U,E)/\Theta(V,E)=\exp(2\pi i\,\alpha(Z_{UV})).
			\]
			\item[(iii)] If $z$ is an extreme point of $P_\Sigma$, then $\alpha(Z_{UV})+\alpha(Z_{VW})=\alpha(Z_{UW})$ for all $U,V,W\in \cI_z(\Sigma)$.
		\end{enumerate}
	\end{lemma}
	
	\begin{proof}
		(i) We first choose $\Omega_i$ as in \eqref{eq:Omega}, and after subtracting even multiples of $[[{T^2}]]$, we may assume their average multiplicities lie in $[0,2]$, and their boundary masses are uniformly bounded. Then by the Poincar\'e inequality and the compactness of functions of bounded variation \cite{AFP00,EvansGariepy15}, there exists a subsequence converging in $L^1$. Moreover, by taking the limit of the subsequence as in \eqref{eq:Omega}, all the limits give the same value of the exponential. This proves continuity.
		
		(ii) Since $\Omega_U\bmod2=\Omega_V\bmod2=E$, $(\Omega_U-\Omega_V)/2$ is an integral current which fills $Z_{UV}$. By the constancy theorem \cite{Federer69}, it differs from $\cA(Z_{UV})$ by an integral multiple of $[[{T^2}]]$. This proves the desired equality.
		
		(iii) Because these three difference cycles belong to $\cG(\mathscr D(U,V,W))$, and by Lemma \ref{lem:iness}(ii), it is inessential. Then the additivity of $\alpha$ on the common inessential set from Proposition \ref{prop:A} shows the equality.
	\end{proof}
	
	The area invariant $\Theta$ depends on both the chosen integral representative $U$ and the filling $E$. We must remove its dependence on $U$ to obtain a map defined on the filling double cover. The following example explains why the signed area relative to a reference cycle is needed.
	
	\begin{example}\label{ex:Theta}
		On $\R^2/\Z^2$, let $E_s=\{s<x<s+1/2\}$, where the first coordinate $x$ is taken modulo $1$. The boundary of $E_s$ consists of two parallel closed geodesics. The path $E_s$, $0\leq s\leq1/2$, joins $E_0$ to $E_0^c$, although $\area(E_s)=1/2$ for all $s$. Thus the unsigned area of a filling does not detect the monodromy of the filling double cover, whereas the signed area relative to a reference cycle, which is essential to the definition of $\Theta$, detects it.
	\end{example}
	
	The following elementary descent lemma is analogous to the vanishing of the first cohomology group for the fine sheaf of continuous real-valued functions; cf. \cite[Chapter II, \S9, Theorems 9.8 and 9.12]{Bredon97}.
	
	\begin{lemma}\label{lem.cocycle}
		Let $p:Y\to C$ be a continuous surjection between compact metric spaces. Suppose that $c$ is a continuous real-valued function on $\{(y,y'):p(y)=p(y')\}$ satisfying $c(y,y')+c(y',y'')=c(y,y'')$ whenever $p(y)=p(y')=p(y'')$. Then there exists a continuous function $h:Y\to\R$ such that
		\[
		c(y,y')=h(y)-h(y')
		\]
		whenever $p(y)=p(y')$.
	\end{lemma}
	
	\begin{proof}[Sketch of proof]
		For each $x\in C$, choose $y_x\in p^{-1}(x)$. By the Tietze extension theorem, the function $c(\,\cdot\,,y_x)$ on $p^{-1}(x)$ extends continuously to $Y$. By continuity, it satisfies the cocycle identity up to an arbitrarily small error on fibers over a neighborhood of $x$. A finite partition of unity on $C$ therefore produces a continuous function $h_1:Y\to\R$ for which $c(y,y')-h_1(y)+h_1(y')$ has arbitrarily small uniform norm on $\{(y,y'):p(y)=p(y')\}$. Iterating with errors tending geometrically to zero gives a uniformly convergent sum $h=\sum_jh_j$ satisfying $c(y,y')=h(y)-h(y')$.
	\end{proof}
	
	Now we are ready to prove (ii) of Proposition \ref{prop_equiv_maps}.
	
	\begin{proof}[Proof of Proposition \ref{prop_equiv_maps}(ii)]
		The common polygon $P_\Sigma$ has a nonzero vertex $z$, which we fix. Then $\cI_z(\Sigma)$ is nonempty for every $\Sigma\in C$. Consider the compact space of integral representatives of class $2z$ and its projection,
		\begin{equation}\label{eq:resolution}
			\{(\Sigma,U):\Sigma\in C,\ U\in \cI_z(\Sigma)\}\to C,\quad (\Sigma,U)\to \Sigma.
		\end{equation}
		We correct $\Theta$ by a continuous function on $\{(\Sigma,U):\Sigma\in C,\ U\in \cI_z(\Sigma)\}$ as follows.
		
		We first claim that the function $(\Sigma,U,V)\to\alpha(Z_{UV})$ is continuous. Indeed, given a convergent sequence $(\Sigma_i,U_i,V_i)\to(\Sigma,U,V)$, let $W_i=(U_i+V_i)/2$ and $Z_i=(U_i-V_i)/2$. The mass measures of $W_i$ and $Z_i$ are mutually singular, and
		\[
		\bM(W_i)+\bM(Z_i)=\bM(\Sigma_i).
		\]
		The same identity holds for the limit. Then by the lower semicontinuity of the mass, we have $\bM(Z_i)\to\bM(Z)$, and Proposition \ref{prop:Acont} proves the claim. Moreover, by Lemma \ref{lem:Theta}(iii), this continuous function satisfies the cocycle identity.
		
		We apply Lemma \ref{lem.cocycle} to the projection \eqref{eq:resolution}, and to the cocycle $c((\Sigma,U),(\Sigma,V))=\alpha(Z_{UV})$, to obtain a continuous function $h$ such that $h(\Sigma,U)-h(\Sigma,V)=\alpha(Z_{UV})$. For $U\in \cI_z(\Sigma)$, we define
		\begin{equation}\label{eq:Psi}
			\Psi(\Sigma,E)=\Theta(U,E)\exp(-2\pi i\,h(\Sigma,U)).
		\end{equation}
		By Lemma \ref{lem:Theta}(ii), this definition does not depend on $U$. The right-hand side is continuous on the compact space of triples $(\Sigma,U,E)$, and the projection of this space onto the filling double cover is a quotient map. Hence $\Psi$ is continuous on the filling double cover. Finally, \eqref{eq:odd} gives $\Psi(\Sigma,E^c)=-\Psi(\Sigma,E)$. This shows that $\Psi$ is the desired map.
	\end{proof}

	\section{The lower bound in the main theorem}\label{S:lower bound}
	
	In this section, we prove the lower bound in the main theorem.
	
	\subsection{Associated polygon}
	
	Recall that, given $\Sigma\in\cZ$, we define the associated polygon in Definition \ref{def:asso polygon}. We discuss several properties of it in this section.
	
	\begin{lemma}\label{lem_per_bound_by_mass}
		For every $\Sigma\in\cZ$ of finite mass, $\Per(P_\Sigma)\leq 2\bM(\Sigma)$.
	\end{lemma}
	
	\begin{proof}
		Fix a unit vector $v\in\R^2$, and choose integral representatives $U_+,U_-$ of $\Sigma$ for which $[U_+]/2$ and $[U_-]/2$ respectively maximize and minimize the pairing with $v$. Write $U_\pm=\sigma_\pm\tau\,\|\Sigma\|$, where $\tau$ is a measurable unit tangent field along $\Sigma$ and $\sigma_\pm\in\{1,-1\}$ almost everywhere. Pairing with the constant $1$-form dual to $v$, and recalling \eqref{eq:W_K},
		\[
		W_{P_\Sigma}(v)=\frac12\int_\Sigma(\sigma_+-\sigma_-)\langle v,\tau\rangle\,ds
		\leq\int_\Sigma|\langle v,\tau\rangle|\,ds.
		\]
		Integrating over unoriented unit vectors and applying Cauchy's formula for the perimeter, we obtain $\Per(P_\Sigma)\leq2\bM(\Sigma)$.
	\end{proof}

	The following theorem gives the obstruction needed for the lower bound.
	
	\begin{theorem}\label{thm_low_lattice_no_sweepout}
		Let $X$ be a finite simplicial complex, and let $\Phi:X\to\cZ$ be an $\mbfF$-continuous map with $\sup_X \bM(\Phi)<\infty$. If
		\begin{equation}\label{eq_PcapLambda}
			\#(P_{\Phi(x)}\cap\Lambda)\leq p\qquad\text{for all }x\in X,
		\end{equation}
		then $(\Phi^*\lambda)^p=0$. In particular, $\Phi$ is not a $p$-sweepout.
	\end{theorem}
	
	We prove the theorem at the end of this section. 
	
	By Lemma \ref{lem_per_bound_by_mass}, $\Phi(X)$ consists of only finitely many classes $[U]/2$, and we denote the collection of these classes by $\scR$. We decompose $\Phi(X)$ into strata, where each stratum consists of slices $\Phi(x)$ with a fixed number of nonzero half-homology classes realized by integral representatives.
	\[
	L_k=\{\Sigma\in \Phi(X):\#\cR(\Sigma)=2k\},\qquad
	\Phi(X)_{\geq k}=\bigcup_{k'\geq k}L_{k'}.
	\]
	We also define $2j=\max_{\Sigma\in \Phi(X)}\#\cR(\Sigma)$ in the following proof. The following elementary observation describes the closure of these strata.

	\begin{lemma}\label{lem:layers}
		For every $S\subset\scR$, the set
		\[
		F_S:=\{\Sigma\in\Phi(X):S\subset\cR(\Sigma)\}
		\]
		is closed. Consequently, $\Phi(X)_{\geq k}$ is closed in $\Phi(X)$, and for every $S\subset\scR$ with $\#S=2k$,
		\[
		E_S:=\{\Sigma\in\Phi(X):\cR(\Sigma)=S\}
		\]
		is closed in $L_k$.
	\end{lemma}
	
	\begin{proof}
		Suppose that $\Sigma_i\to\Sigma$ in the $\mbfF$-topology and $S\subset\cR(\Sigma_i)$. For each $z\in S$, choose $U_i\in\cI_z(\Sigma_i)$. By Lemma \ref{lem:closedgraph}, after passing to a subsequence, $U_i$ converges to an element of $\cI_z(\Sigma)$. Hence $S\subset\cR(\Sigma)$, and $F_S$ is closed.
		
		Since $\scR$ is finite,
		\[
		\Phi(X)_{\geq k}
		=
		\bigcup_{\substack{S\subset\scR\\ \#S=2k}}F_S
		\]
		is closed. If $\#S=2k$, then $E_S=F_S\cap L_k$, which is closed.
	\end{proof}

	Now we are ready to prove Theorem \ref{thm_low_lattice_no_sweepout}. The proof is to construct an equivariant map factored through the double cover.
	
	\begin{proof}[Proof of Theorem \ref{thm_low_lattice_no_sweepout}]
		In the proof, we use the widetilde to denote the pullback of the filling double cover. We let $Y_0=S^0\cong\{\pm 1\}$ and we let $Y_k=S^1$ for $k\geq 1$. We construct open sets $W_k\subset\Phi(X)$ and equivariant maps from the filling double cover over $W_k$ to $Y_k$ such that
		\[
		\bigcup_{l\geq k}W_l\supset\Phi(X)_{\geq k},
		\]
		by backward induction on $k$.
		
		Suppose we have constructed $W_{k+1},\cdots,W_j$. Since they cover $\Phi(X)_{\geq k+1}$, the set
		\[
		C_k=\Phi(X)_{\geq k}\backslash\bigcup_{l>k}W_l
		\]
		is compact and contained in $L_k$. By Lemma \ref{lem:layers}, the finitely many nonempty sets $C_k\cap E_S$, $\#S=2k$, are compact and pairwise disjoint. For each $S$ with $\#S=2k$ and $C_k\cap E_S\neq\emptyset$, Proposition \ref{prop_equiv_maps} gives an equivariant map
		\[
		f_{k,S}:\widetilde{C_k\cap E_S}\to Y_k.
		\]
		By the Tietze extension theorem and odd symmetrization, after shrinking to pairwise disjoint neighborhoods $W_{k,S}$ of $C_k\cap E_S$, each $f_{k,S}$ extends to an equivariant map $\widetilde W_{k,S}\to Y_k$. Let us define
		\[
		W_k=\bigcup_S W_{k,S},
		\]
		and we define an equivariant map combining all these maps
		\[
		f_k:\widetilde W_k\to Y_k.
		\]
		Thus the backward induction gives
		\[
		\bigcup_{l\geq k}W_l\supset\Phi(X)_{\geq k},
		\]
		and in particular $W_0,\cdots,W_j$ cover $\Phi(X)$.
		
		Let us choose a partition of unity $\rho_0,\cdots,\rho_j$ on $X$ subordinate to $\{\Phi^{-1}(W_k)\}$. Let $\widetilde X$ be the pullback of the filling double cover, and denote the induced local maps again by $f_k$. Viewing $f_0$ in $S^0\subset\R$ and $f_k$ in $S^1\subset\R^2$ for $k\geq1$, define
		\begin{equation}\label{eq:join}
			F(\tilde x)
			=
			\left(
			\sqrt{\rho_0(x)}f_0(\tilde x),
			\cdots,
			\sqrt{\rho_j(x)}f_j(\tilde x)
			\right)
			\in S^{2j}.
		\end{equation}
		This map is continuous and equivariant with respect to the antipodal actions. Hence it descends to a map $X\to \mathbb{RP}^{2j}$ whose pullback of the antipodal double cover is the filling double cover. Therefore
		\[
		(\Phi^*\lambda)^{2j+1}=0.
		\]
		
		Finally, let us choose $\Sigma\in\Phi(X)$ with $\#\cR(\Sigma)=2j$. By \eqref{eq:latticecount} and \eqref{eq_PcapLambda},
		\[
		2j+1\leq\#(P_\Sigma\cap\Lambda)\leq p.
		\]
		Thus $(\Phi^*\lambda)^p=0$.
	\end{proof}
	
	\subsection{The lower bound}
	Finally, we prove the lower bound for $\omega_p(\Lambda)$.
	
	\begin{proof}[Proof of $\geq$ in Theorem \ref{main_thm_omega_p}]
		Let the right-hand side of \eqref{eq:main} be $C_p(\Lambda)$. Suppose that $\omega_p(\Lambda)<C_p(\Lambda)$. Then there exists an $\mbfF$-continuous $p$-sweepout $\Phi:X\to\cZ$ with $\sup_X\bM(\Phi)<C_p(\Lambda)$ \cite{MarquesNeves17,LMN}. By Lemma \ref{lem_per_bound_by_mass}, for every $x\in X$,
		\[
		\Per(P_{\Phi(x)})\leq 2\bM(\Phi(x))<2C_p(\Lambda).
		\]
		Since $P_{\Phi(x)}$ is a centrally symmetric convex lattice polygon centered at the origin, Proposition \ref{prop_C_p_bound_by_perimeter} implies that $\#(P_{\Phi(x)}\cap\Lambda)\leq p$. Then Theorem \ref{thm_low_lattice_no_sweepout} shows that $\Phi$ is not a $p$-sweepout, which is a contradiction.
	\end{proof}
	
	\section{Widths of finite isometric quotients of flat tori}\label{S:widths of quotient}
	In this section, we prove Theorem \ref{thm:affine}. The arguments are essentially the same as in Sections \ref{S:upper bound}--\ref{S:lower bound} after taking the group action into account. Therefore, we only highlight the additional ingredients needed for the quotient setting.
	
	\subsection{Preparation for the quotient formula}
	
	Throughout, we fix a lattice $\Lambda$ in $\R^2$, $T^2=\R^2/\Lambda$, an isometric group $G$ of $T^2$, and $Q=T^2/G$. For every element $g\in G$, we can decompose it $g(x)=L_g(x)+t_g$, where $L_g\in O(2)$ preserves $\Lambda$ and $t_g$ is the translation defined modulo $\Lambda$. We denote by $\bar G$ the group of linear parts. We will consider relative cycles along the boundary of $Q$ and we view the cone points as interior points.
	
	Let $\Lambda^*=\{\xi\in\R^2:\langle\xi,z\rangle\in\Z\text{ for all }z\in\Lambda\}$ be the dual lattice. Since $\det(z,x)=\langle Jz,x\rangle$, we have $\Lambda^*=J\Lambda/A$, so the Fourier modes below are those we used in Section \ref{S:upper bound}. For $\xi\in\Lambda^*$, the function $e_\xi(x)=e^{2\pi i\langle\xi,x\rangle}$ is well defined on $T^2$, and for $g\in G$ we have
	\begin{equation}\label{eq:plane-wave}
		e_\xi\circ g=e^{2\pi i\langle\xi,t_g\rangle}e_{L_g^{-1}\xi}.
	\end{equation}
	
	The computation of the widths relies on a compatible condition on the lattice.
	
	\begin{definition}\label{def:compatible}
		A vector $\xi\in\Lambda^*$ is called $G$-compatible if
		\[
		\langle\xi,t_g\rangle\in\Z
		\]
		whenever $g\in G$ satisfies $L_g\xi=\xi$. We denote the set of $G$-compatible vectors by $\Lambda_G^*$.
	\end{definition}
	
	Let $\cK_G$ be the set of centrally symmetric, $\bar G$-invariant convex polygons with vertices in $\Lambda^*$, and we include segments and $\{0\}$ as degenerate polygons. For $K\in\cK_G$, We define
	\[
	V_K=\Bigl\{\sum_{\xi\in K\cap\Lambda^*}a_\xi e_\xi:\ a_{-\xi}=\overline{a_\xi}\Bigr\}
	\]
	be the real vector space of trigonometric polynomials with frequencies in $K$, as in Section \ref{S:upper bound} for the no quotient case. Let $V_K^G$ be the subspace of $G$-invariant functions in $V_K$, and we define $q_G(K)=\dim_\R V_K^G$. Since $J\Lambda_G/A=\Lambda_G^*$, Lemma \ref{lem:fourier-count}(ii) below shows that $q_G(JP/A)$ is the same as $q_G(P)$ for every $G$-invariant lattice polygon $P$.
	
	The next lemma shows that we can restrict our attention to the set of $G$-compatible frequencies in Theorem \ref{thm:affine}.
	
	\begin{lemma}\label{lem:fourier-count}
		The following statements hold for every $K\in\cK_G$. 
		\begin{enumerate}[label=(\roman*)]
			\item The set $\Lambda_G^*$ is invariant under $\bar G$ and under $\xi\to-\xi$.
			\item $q_G(K)
			=
			\#((K\cap\Lambda_G^*)/\bar G)$. In particular, if $G$ fixes the origin, then $\Lambda_G^*=\Lambda^*$.
			\item
			\begin{equation}\label{eq:trace}
				q_G(K)
				=
				\frac1{|G|}
				\sum_{g\in G}
				\sum_{\substack{\xi\in K\cap\Lambda^*\\L_g\xi=\xi}}
				e^{2\pi i\langle\xi,t_g\rangle}.
			\end{equation}
			\item If $K'$ is the convex hull of $K\cap\Lambda_G^*$, then
			\[
			K'\in\cK_G,\qquad
			q_G(K')=q_G(K),\qquad
			\Per(K')\leq\Per(K).
			\]
			In particular, in Theorem \ref{thm:affine}, the minimum may be taken over polygons whose vertices are $G$-compatible.
		\end{enumerate}
	\end{lemma}
	
	\begin{proof}
		We first prove (iii). Let $W_K$ be the complex span of $\{e_\xi:\xi\in K\cap\Lambda^*\}$, and by \eqref{eq:plane-wave}, the action of $g$, denoted by $R_g$, on the Fourier modes preserves $W_K$. Because the $g$-action commutes with the complex conjugate, $q_G(K)=\dim_\R V_K^G=\dim_{\mathbb C}W_K^G$. Thus, the averaging operator $\Pi=\frac1{|G|}\sum_{g\in G}R_g$ is a projection of $W_K$ onto $W_K^G$, and hence,
		\[
		q_G(K)=\tr\Pi=\frac1{|G|}\sum_{g\in G}\tr R_g.
		\]
		Suppose we write $\Pi$ in the basis $\{e_\xi\}$, then the diagonal coefficient of $R_g$ at $e_\xi$ is $e^{2\pi i\langle\xi,t_g\rangle}$ if $L_g\xi=\xi$, and zero otherwise. This proves (iii).
		
		Next we prove (i) and (ii). For $\xi\in\Lambda^*$, let $G_\xi=\{g\in G:L_g\xi=\xi\}$ be the stabilizer of $\xi$, and $S_\xi=\sum_{g\in G_\xi}e^{2\pi i\langle\xi,t_g\rangle}$. If $g,h\in G_\xi$, then $t_{hg}=L_ht_g+t_h$ modulo $\Lambda$ and $\langle\xi,L_ht_g\rangle=\langle\xi,t_g\rangle$, so
		\[
		e^{2\pi i\langle\xi,t_{hg}\rangle}=e^{2\pi i\langle\xi,t_h\rangle}e^{2\pi i\langle\xi,t_g\rangle}.
		\]
		If $\xi\in\Lambda_G^*$, every term of $S_\xi$ equals $1$, so $S_\xi=|G_\xi|$. Otherwise, we choose $h\in G_\xi$ with $e^{2\pi i\langle\xi,t_h\rangle}\neq1$. Since $g\to hg$ permutes $G_\xi$, the identity above gives $S_\xi=e^{2\pi i\langle\xi,t_h\rangle}S_\xi$, and hence $S_\xi=0$. Moreover, for $h\in G$, the map $g\to hgh^{-1}$ sends $G_\xi$ onto $G_{L_h\xi}$, and $\langle L_h\xi,t_{hgh^{-1}}\rangle\equiv\langle\xi,t_g\rangle\pmod\Z$. Thus $L_h\xi$ is compatible whenever $\xi$ is. For $\xi\to-\xi$, the invariance is straightforward. Thus (i) holds. Then we sum \eqref{eq:trace} over $g\in G_\xi$ to obtain
		\[
		q_G(K)=\frac1{|G|}\sum_{\xi\in K\cap\Lambda^*}S_\xi=\sum_{\xi\in K\cap\Lambda_G^*}\frac{|G_\xi|}{|G|}=\#\bigl((K\cap\Lambda_G^*)/\bar G\bigr).
		\]
		The last equality holds because each point $L_g\xi$ of the orbit $\bar G\xi$ is obtained from exactly $|G_\xi|$ elements of $G$, and $K\cap\Lambda_G^*$ is a union of such orbits by (i). This proves (ii).
		
		Finally, (iv) follows from (i), (ii), and the monotonicity of perimeter under inclusion.
	\end{proof}
	
	\subsection{Lifting sweepouts}
	
	Given a sweepout of $Q$, we need to lift it to a $G$-invariant family on $T^2$. Let $\pi_Q:T^2\to Q$ be the quotient map. Pullback by $\pi_Q$ sends finite-perimeter sets in $Q$ to $G$-invariant finite-perimeter sets in $T^2$, up to sets with measure $0$. Moreover, the pullback commutes with taking the complement. 
	
	For every finite-perimeter set $E\subset Q$ and every open set $U\subset Q$, we have
	\begin{equation}\label{eq:lifting}
		|\pi_Q^{-1}(E)|=|G|\,|E|,
		\qquad
		\Per(\pi_Q^{-1}(E);\pi_Q^{-1}(U))=|G|\,\Per(E;U\backslash\partial Q).
	\end{equation}
	Indeed, $\pi_Q$ is a covering of degree $|G|$ away from the boundaries and cone points. By reflection invariance, $\pi_Q^{-1}(E)$ has the same boundary values on both sides of each boundary, so its perimeter does not count the mirrors. Also, the cone points are isolated and do not contribute to the perimeter.
	
	For families of bounded mass, \eqref{eq:lifting} and the compactness theorem for functions of bounded variation \cite{AFP00,EvansGariepy15} show that lifting preserves the filling double covers, the flat continuity, and nonconcentration of mass. We can also apply \eqref{eq:lifting} to symmetric differences to see that lifting preserves mass continuity.

	\eqref{eq:lifting} also addresses the following technical point. The lower bound in Theorem \ref{thm:affine} requires the sweepouts to be continuous in mass, not only in the flat topology. If the ambient manifold is smooth with possibly nonempty smooth boundary, we apply the relative interpolation theorem \cite{MarquesNeves17,LMN}, and since sweepouts have no concentration of mass, the resulting increase in maximal mass tends to zero with the size of the neighborhoods. Thus, up to a very tiny change of maximal mass, we may assume that a sweepout of $Q$ is mass-continuous. For $Q$ which possibly contains cone points and corners, we smooth the finitely many cone points and corners of $Q$ in arbitrarily small neighborhoods. By \eqref{eq:lifting}, its lift is a mass-continuous $G$-invariant family on $T^2$, with all masses multiplied by $|G|$.

	\subsection{Proof of the quotient formula}
	
	Now we adapt the idea of the proof in Section \ref{S:lower bound} to cycles with $G$-invariant fillings. We need two new ingredients: the action of $G$ on the area invariant, which forces the vertices of the associated polygon to be compatible; and a local map to $S^0$ when $\bar G$ reverses a vertex.
	
	Suppose that $\Sigma\in\cZ$ has a $G$-invariant filling $E$. Then the action of each group element, $g_\#$, preserves $\cI(\Sigma)$ and satisfies $[g_\#U]=L_g[U]$. Hence, $P_\Sigma$ is $\bar G$-invariant, and as a consequence, the polygon $JP_\Sigma/A$ is also $\bar G$-invariant. Also, if we lift $g$ to a map in $\R^2$, it maps the bounded Jordan fillings to bounded Jordan fillings. Therefore, Proposition \ref{prop:A} gives
	\begin{equation}\label{quot:eq:alpha-natural}
		\cA(g_\#Z)=g_\#\cA(Z),
		\qquad \alpha(g_\#Z)=(\det L_g)\alpha(Z)
	\end{equation}
	for every integral cycle $Z$ whose mass measure is concentrated on an inessential set. 
	
	Next we try to define the area invariant under $G$-symmetry. We need to understand how the vertices of the associated polygon change under the group action. For a nonzero vertex $z$ of $P_\Sigma$, let $\xi_z=Jz/A$, and we choose $R_z=\pi_\#[[0,z]]$. Let $z$ be a nonzero vertex of $P_\Sigma$, and suppose that $L_g\xi_z=s\xi_z$ with $s\in\{1,-1\}$. Since $L_gJ=(\det L_g)JL_g$, this means $L_gz=s(\det L_g)z$. For $U\in\cI_z(\Sigma)$, the cycle $g_\#U$ is again an integral representative of $\Sigma$, but its class is $2L_gz=2s(\det L_g)z$. Therefore, we need to define
	\begin{equation}\label{eq:T_gU}
		T_gU:=s(\det L_g)\,g_\#U,
	\end{equation}
	so that $[T_gU]=2z$ and $T_gU\in\cI_z(\Sigma)$. Thus $T_g$ is the action of $g$ on integral representatives, with possibly a sign flip that keeps the homology class fixed.

	\begin{lemma}\label{lem:G-area_invariant}
		Let $z$ be a nonzero vertex of $P_\Sigma$, and suppose that $L_g\xi_z=s\xi_z$ for some $g\in G$ and $s\in\{1,-1\}$. Then for every $U\in\cI_z(\Sigma)$,
		\begin{equation}\label{eq:G-area_invariant}
			\Theta(T_gU,E)=e^{-2\pi i\langle\xi_z,t_g\rangle}\Theta(U,E)^s.
		\end{equation}
		Moreover, every vertex of $JP_\Sigma/A$ is $G$-compatible.
	\end{lemma}

	\begin{proof}
		We consider the projected parallelogram $B$ parametrized by $(r,u)\to rt_g+uz$ on $[0,1]^2$. Then it satisfies
		\[
		\partial B=(s\det L_g) g_\#R_z-R_z,
		\qquad B(dA)=-A\langle\xi_z,t_g\rangle.
		\]
		If $\Omega$ defines $\Theta(U,E)$, then $(s\det L_g) g_\#\Omega+2B$ has boundary $T_gU-2R_z$ and equals $g(E)=E$ modulo $2$. Since $(g_\#\Omega)(dA)=(\det L_g)\Omega(dA)$, and hence $((s\det L_g)g_\#\Omega)(dA)=s\,\Omega(dA)$,
		we evaluate its signed area to get \eqref{eq:G-area_invariant}.
		
		To show that $JP_\Sigma/A$ is $G$-compatible, let $s=1$ and let $N$ be the order of $g$, so that $T_g^N=1$. By Lemma \ref{lem:iness}(ii), the differences of $U,T_gU,\cdots,T_g^{N-1}U$ have mass measures concentrated on a common inessential set. By \eqref{quot:eq:alpha-natural}, the numbers $\alpha(Z_{T_g^jU,T_g^{j+1}U})$ are all equal, and by Lemma \ref{lem:Theta}(iii) their sum is $\alpha(Z_{U,U})=0$. Thus $\alpha(Z_{U,T_gU})=0$, and Lemma \ref{lem:Theta}(ii) together with \eqref{eq:G-area_invariant} gives $\langle\xi_z,t_g\rangle\in\Z$, which implies that $\xi_z$ is $G$-compatible.
	\end{proof}
	
	Let $\scP_G$ be the set of centrally symmetric, $\bar G$-invariant convex lattice polygons $P$ centered at the origin, including the line segments and $\{0\}$, for which every vertex of $JP/A$ is $G$-compatible. By Lemma \ref{lem:G-area_invariant}, $P_\Sigma\in\scP_G$ whenever $\Sigma$ has $G$-invariant fillings. For $P\in\scP_G$, we write  $q(P)=q_G(JP/A)$, and define
	\begin{equation}\label{quot:eq:PG}
		d(P)=
		\begin{cases}
			1,&\text{if }P=\{0\}\text{, or if }-\xi\in\bar G\xi\text{ for some nonzero vertex }\xi\text{ of }JP/A,\\
			2,&\text{otherwise}.
		\end{cases}
	\end{equation}
	Note that $d(P)=2$ exactly when $P\neq\{0\}$ and no element of $\bar G$ maps a nonzero vertex $\xi$ of $JP/A$ to $-\xi$. This number will determine whether the local equivariant map has target space $S^0$ or $S^1$. 
	
	In fact, the case $-\xi\in\bar G\xi$ is the only one that changes the local construction. If $\xi=Jz/A$, $g$ maps the family of closed geodesics parallel to $z$ to itself but reverses its transverse orientation as a reflection. Then $\xi$ and $-\xi$ lie in the same $\bar G$-orbit, so such a vertex adds only one orbit to the count of $q$, rather than the two orbits of $\xi$ and $-\xi$. In this case, however, \eqref{eq:G-area_invariant} with $s=-1$ relates $\Theta$ to its inverse, which is its complex conjugate. Then the circle-valued phase of Proposition \ref{prop_equiv_maps}(ii) can be corrected to take values in $S^1\cap\R=S^0$.
	
	\begin{lemma}\label{lem:local}
		Let $C\subset\cZ$ be an $\mbfF$-compact family of bounded mass whose cycles have $G$-invariant fillings and a common associated polygon $P$. Then its filling double cover admits a continuous equivariant map to $S^{d(P)-1}$.
	\end{lemma}
	
	\begin{proof}
		The case $P=\{0\}$ is Proposition \ref{prop_equiv_maps}(i). If $d(P)=2$, the proof of Proposition \ref{prop_equiv_maps}(ii) applies, since it uses only a common nonzero vertex. Thus, the main challenge is the other possibility of $d(P)=1$.
		
		In this case, there exists a vertex $z$ of $P$ and $g\in G$ with $L_g\xi_z=-\xi_z$. Let us define $h(\Sigma,U)=\frac12\alpha(Z_{U,TU})$. Its continuity follows from the proof of Proposition \ref{prop_equiv_maps}(ii). Since $T_g=-(\det L_g)g_\#$, \eqref{quot:eq:alpha-natural} shows that $\alpha(Z_{TU,TV})=-\alpha(Z_{UV})$, and then, Lemma \ref{lem:Theta}(iii) shows that $h(\Sigma,U)-h(\Sigma,V)=\alpha(Z_{UV})$. Hence, as in \eqref{eq:Psi},
		\begin{equation}\label{quot:eq:real-phase}
			\Psi(\Sigma,E)=e^{\pi i\langle\xi_z,t_g\rangle}\Theta(U,E)e^{-2\pi i h(\Sigma,U)}
		\end{equation}
		is independent of $U$, continuous, and odd under taking the complement. Finally, Lemma \ref{lem:Theta}(ii) and \eqref{eq:G-area_invariant} implies that $\Theta(U,E)^2=e^{4\pi i h(\Sigma,U)-2\pi i\langle\xi_z,t_g\rangle}$, so $\Psi^2=1$, which implies that $\Psi$ maps onto $S^0$.
	\end{proof}
	
	Now we state the analog of Theorem \ref{thm_low_lattice_no_sweepout}. The vertex count in Section \ref{S:lower bound}, especially in Theorem \ref{thm_low_lattice_no_sweepout}, should be understood as follows. If $P,P'\in\scP_G$ and $P'\subsetneqq P$, then
	\begin{equation}\label{eq:orbit-drop}
		q(P)-q(P')\geq d(P).
	\end{equation}
	In fact, by Lemma \ref{lem:fourier-count}(ii), $q$ counts $\bar G$-orbits of compatible frequencies. If some vertex $\xi$ of $JP/A$ lies outside $JP'/A$, by Lemma \ref{lem:fourier-count}(i), its orbit contains $G$-compatible vertices outside $JP'/A$. When $d(P)=2$, the orbit of $-\xi$ is a second such orbit. This explains why we count as in \eqref{eq:orbit-drop}.
	
	\begin{theorem}\label{thm:obstruction}
		Let $\Phi:X\to\cZ$ be an $\mbfF$-continuous family of bounded mass cycles. If its cycles have $G$-invariant fillings and $q(P_{\Phi(x)})\leq p$ for every $x\in X$, then $(\Phi^*\lambda)^p=0$.
	\end{theorem}
	
	\begin{proof}
		By Lemma \ref{lem_per_bound_by_mass}, $\{P_{\Phi(x)}\}_{x\in X}$ consists of only finitely many polygons. Then we use the backward induction construction as in Theorem \ref{thm_low_lattice_no_sweepout}, but we index the strata by these polygons and order them by decreasing $q(P)$. Just as in the proof of Theorem \ref{thm_low_lattice_no_sweepout}, the Tietze extension theorem and odd symmetrization extend each map given by Lemma \ref{lem:local} from the filling double cover over a compact stratum to the filling double cover over an open neighborhood. This gives an open cover $\{W_P\}$ of $\Phi(X)$ and equivariant maps $f_P:\widetilde W_P\to S^{d(P)-1}$. We also require
		\[
		\overline W_P\subset\{\Sigma\in\Phi(X):P_\Sigma\subset P\},
		\qquad
		\overline W_P\cap\overline W_{P'}=\emptyset
		\quad\text{if }P,P'\text{ are incomparable}.
		\]

		For each polygon $P$ that occurs, let $\iota_P:\R^{d(P)}\to\R^p$ be the inclusion sending
		\[
		(x_{q(P)-d(P)+1},\cdots,x_{q(P)})\to (0,0,\cdots,0,x_{q(P)-d(P)+1},\cdots,x_{q(P)},0,\cdots,0).
		\]
		The left-hand side is a block of a single coordinate, or two coordinates. If $W_P\cap W_{P'}\neq\emptyset$, then $P$ and $P'$ are comparable, say $P'\subsetneq P$, and \eqref{eq:orbit-drop} implies that $q(P')\le q(P)-d(P)$. Thus, the label collections $\{q(P)-d(P)+1,\cdots,q(P)\}\cap \{q(P')-d(P')+1,\cdots,q(P')\}=\emptyset$. Let $\{\rho_P\}$ be a partition of unity on $X$ subordinate to $\{\Phi^{-1}(W_P)\}$, denote the pulled-back maps on $\widetilde X$ again by $f_P$, and we define
		\[
		F(\tilde x)=\sum_P\sqrt{\rho_P(x)}\,\iota_P\bigl(f_P(\tilde x)\bigr),
		\]
		where $x$ is the image of $\tilde x$ in $X$ and a term is zero when $\rho_P(x)=0$. Then $|F(\tilde x)|^2=\sum_P\rho_P(x)=1$, which implies that $F:\widetilde X\to S^{p-1}$ is continuous and equivariant. It induces a map $X\to\mathbb{RP}^{p-1}$ that pulls back the antipodal double cover to the filling double cover, so $\Phi^*\lambda$ is pulled back from $H^1(\mathbb{RP}^{p-1};\Z_2)$, and hence, $(\Phi^*\lambda)^p=0$.
	\end{proof}
	
	\begin{proof}[Proof of Theorem \ref{thm:affine}]
		By Lemma \ref{lem:fourier-count}(iv), we may restrict to polygons with compatible vertices.
		
		For the upper bound, apply the proof of Theorem \ref{thm:torus-upper} to $V_K^G$. Its nodal family projects to a $(q_G(K)-1)$-sweepout of $Q$, and \eqref{eq:fourier-mass} and \eqref{eq:lifting} give maximal mass at most $\area(Q)\Per(K)/2$.
		
		For the lower bound, let $W$ denote the right-hand side of \eqref{quot:eq:main}. Suppose by contradiction that $\omega_p(Q)<W$. Then, by the interpolation and lifting arguments above, there exists a $p$-sweepout of $Q$ with maximal mass less than $W$, and it gives a mass-continuous lift $\Phi$ with invariant fillings, $(\Phi^*\lambda)^p\neq0$, and $\sup_X\bM(\Phi)<|G|W$. Then Lemmas \ref{lem:G-area_invariant} and \ref{lem_per_bound_by_mass} imply
		\[
		JP_{\Phi(x)}/A\in\cK_G,
		\qquad
		\Per(JP_{\Phi(x)}/A)\leq\frac{2\bM(\Phi(x))}{A}<\frac{2|G|W}{A}.
		\]
		But $W$ is the right-hand side of \eqref{quot:eq:main}, so we have $q(P_{\Phi(x)})\leq p$, which is a contradiction to Theorem \ref{thm:obstruction}.
	\end{proof}

	\appendix
	
	\section{Application: the widths determine a flat torus}\label{S:widths determine a flat torus}
	We now show that the sequence of $p$-widths determines a flat two-dimensional torus up to isometry.
	\begin{theorem}\label{thm:flat_rigidity}
		Two flat two-dimensional tori with the same volume spectrum are isometric.
	\end{theorem}
	
	\subsection{Lattice polygons and normalization}
	To prove Theorem \ref{thm:flat_rigidity}, it would be helpful to use some language in the study of convex lattice polygons. Let $\cE$ be the set of unoriented primitive lattice directions. We will fix an ordering of $\cE$ and a primitive representative $e$ for each direction.
	
	A family $r=\{r_e\}_{e\in\cE}$ of nonnegative integers is called {\bf admissible} if it has finitely many nonzero members and
	\[
	\sum_e r_e e\in2\Lambda.
	\]
	For such a family, define the {\bf zonotope}
	\[
	Z(r)=\sum_e\left[-\frac{r_e}{2}e,\frac{r_e}{2}e\right]
	\]
	and the quantities
	\[
	M(r)=\sum_e r_e|e|,
	\qquad
	D(r)=\sum_e r_e+\frac1A\sum_{e<f}r_er_f|\det(e,f)|.
	\]
	We will show that $M(r)$ is half the perimeter of the zonotope $Z(r)$, and $D(r)+1$ is the number of lattice points in $Z(r)$. The admissible condition and these definitions are independent of the choice of signs of the representatives.
	
	In the following, we use the language of convex geometry to give an alternative characterization of the $p$-width.
	
	\begin{definition}\label{def.C_p}
		For each integer $p\geq1$, let
		\[
		C_p(\Lambda)=\min\{M(r):r\text{ is admissible and }D(r)\geq p\}.
		\]
	\end{definition}
	
	\begin{proposition}\label{prop_C_p_bound_by_perimeter}
		Every centrally symmetric convex lattice polygon centered at the origin is $Z(r)$ for an admissible family $r$. Conversely, for every admissible family $r$, the zonotope $Z(r)$ is a centrally symmetric convex lattice polygon and satisfies
		\[
		\frac12\Per(Z(r))=M(r),
		\qquad
		\#(Z(r)\cap\Lambda)=D(r)+1.
		\]
		Consequently, $C_p(\Lambda)$ equals half the minimum perimeter among such polygons containing at least $p+1$ lattice points.
	\end{proposition}
	
	\begin{proof}
		Let $r$ be admissible. Every vertex of $Z(r)$ is a signed half-sum of the vectors $r_e e$, and hence belongs to $\Lambda$ by the condition $\sum_e r_e e\in2\Lambda$. If $Z(r)$ has nonempty interior, then for each $e$ with $r_e>0$, the two edges parallel to $e$ both have lattice length $r_e$. Thus
		\[
		\Per(Z(r))=2\sum_e r_e|e|,
		\qquad
		\#(\partial Z(r)\cap\Lambda)=2\sum_e r_e.
		\]
		The area formula for zonotopes \cite{Shephard74,Schneider14} shows that $\area(Z(r))=\sum_{e<f}r_er_f|\det(e,f)|$, and Pick's theorem \cite{Pick1899} shows that
		\begin{equation}\label{eq:zonotope-count}
			\#(Z(r)\cap\Lambda)
			=1+\frac{\area(Z(r))}{A}+\sum_e r_e
			=D(r)+1.
		\end{equation}

		Conversely, let $P$ be a centrally symmetric convex lattice polygon with nonempty interior, centered at the origin. Let us consider, instead of the vertices, the edges. Since $P$ is centrally symmetric, its oriented edge vectors can be written, in cyclic order, as $r_{e_1}e_1,\cdots,r_{e_m}e_m,
		-r_{e_1}e_1,\cdots,-r_{e_m}e_m$. Then the zonotope $
		Z(r)=\sum_e\left[-\frac{r_e}{2}e,\frac{r_e}{2}e\right]$ has exactly the same oriented edge vectors, in the same cyclic order. Thus, $P$ and $Z(r)$ are the same. Because each vertex is a signed half-sum of these edge vectors, we have $\sum_e r_e e\in2\Lambda$. Thus $r$ is admissible.
	\end{proof}
	
	Together with Theorem \ref{main_thm_omega_p}, Proposition \ref{prop_C_p_bound_by_perimeter} shows that $C_p(\Lambda)=\omega_p(\Lambda)$ for every $p\geq1$. 
	
	Let us compute the first two widths explicitly.
	\begin{lemma}\label{lem:1stw_sys}
		We have $\omega_1(\Lambda)=\omega_2(\Lambda)=2\sys(\Lambda)$.
	\end{lemma}
	
	\begin{proof}
		Because a nontrivial centrally symmetric lattice polygon must contain a segment $[-z,z]$ with $z\in\Lambda\backslash\{0\}$, it must have perimeter at least $4|z|\geq4\sys(\Lambda)$. If $z$ is a shortest nonzero lattice vector, then $[-z,z]$ contains three lattice points and has perimeter $4\sys(\Lambda)$. By Theorem \ref{main_thm_omega_p}, we get the lemma.
	\end{proof}
	
	Therefore, the first width determines the systole. Since the widths scale linearly under homotheties, we may normalize the systole to $1$. After a rotation, and possibly a reflection, we may write a basis for the lattice
	\[
	u=(1,0),\qquad v=(x,A),\qquad 0\leq x\leq\frac12,
	\]
	where $A$ is the area of the torus. Equation \eqref{eq:Weyl} gives $A=\lim_{p\to\infty}\frac{C_p(\Lambda)^2}{\pi p}$, so the spectrum determines $A$. Since $|v|\geq1$, we also have $x^2+A^2\geq1$ and $A\geq\sqrt3/2$. Thus, to show the widths determine the flat metric on $T^2$, we only need to show that the widths determine $x$.
	
	For fixed $A$, the lengths of the vectors $v+tu$ are given by
	\begin{equation}\label{eq:lenrt}
		|v+tu|=\sqrt{(x+t)^2+A^2},
		\qquad t\in\R.
	\end{equation}
	Later, we will use this quantity to study the polygon that realizes the widths. We will use the identity
	\begin{equation}\label{eq.deriv.|v+tu|}
		\frac{\partial}{\partial t}|v+tu|
		=\frac{\partial}{\partial x}|v+tu|
		=\psi(x+t),
		\qquad
		\psi(s)=\frac{s}{\sqrt{s^2+A^2}},
	\end{equation}
	where $\psi$ is odd and strictly increasing. For integers $k\geq l$, we have
	\[
	\frac{\det(v+ku,v+lu)}{A}=k-l.
	\]
	The next lemma minimizes the sum of the lengths for a fixed value of this determinant.
	
	\begin{lemma}\label{lem:bal}
		Let $k\geq l$ be integers with $k-l=n$. Then
		\[
		|v+ku|+|v+lu|
		\geq
		\begin{cases}
			|v+mu|+|v-mu|,& n=2m,\\
			|v+mu|+|v-(m+1)u|,& n=2m+1.
		\end{cases}
		\]
		Moreover, $|v+mu|+|v-mu|$ is nondecreasing in $m\geq0$.
	\end{lemma}
	
	\begin{proof}
		Write $k=c+n/2$ and $l=c-n/2$, where $c\in n/2+\Z$. By \eqref{eq:lenrt}, the sum $|v+(c+n/2)u|+|v+(c-n/2)u|$ is strictly convex in $c$ and symmetric about $c=-x$. Its minimum over $n/2+\Z$ is attained at a point closest to $-x$. Since $0\leq x\leq1/2$, we may choose $c=0$ when $n$ is even and $c=-1/2$ when $n$ is odd. This proves the two inequalities. Finally, \eqref{eq.deriv.|v+tu|} gives
		\[
		\frac{d}{dm}\bigl(|v+mu|+|v-mu|\bigr)
		=\psi(x+m)-\psi(x-m)\geq0
		\]
		for $m\geq0$.
	\end{proof}

	\subsection{Polygons with five lattice points}
	
	Let
	\[
	\Delta:=|v-u|+|v+u|
	\]
	be the sum of the lengths of the two diagonals of the fundamental parallelogram.
	
	\begin{lemma}\label{lem:key}
		If $x\in[0,1/2]$ and $\sqrt{x^2+A^2}\geq 1$, then $\Delta\leq 1+2A$. Consequently,
		\begin{equation}\label{eq:key}
			\Delta\leq4A,\qquad \Delta\leq1+|2v-u|.
		\end{equation}
	\end{lemma}
	
	\begin{proof}
		Since $0\leq x\leq\frac12$, we have $|v\pm u|>A$, and therefore
		\[
		\partial_A(\Delta-2A)=\frac{A}{|v-u|}+\frac{A}{|v+u|}-2<0.
		\]
		It therefore suffices to consider the least admissible value $A=\sqrt{1-x^2}$. In this case, $|v\mp u|=\sqrt{2\mp2x}$, so
		\[
		\Delta-2A=2\sqrt{1+\sqrt{1-x^2}}-2\sqrt{1-x^2}.
		\]
		The function $t\to2\sqrt{1+t}-2t$ is decreasing for $t\geq0$ and equals $1$ at $t=\sqrt3/2$. Since $\sqrt{1-x^2}\geq\sqrt3/2$, we obtain $\Delta-2A\leq1$. The two inequalities in \eqref{eq:key} follow from $2A\geq\sqrt3>1$ and from $|2v-u|\geq 2A$, respectively.
	\end{proof}
	
	\begin{proposition}\label{prop:C3}
		$C_3=C_4=\min\{4,\,\Delta\}$. The value $4$ is realized by the segment $[-2u,2u]$, and the value $\Delta$ is realized by the parallelogram with vertices $\pm u,\pm v$.
	\end{proposition}
	
	\begin{proof}
		The equality $C_3=C_4$ follows from Corollary \ref{cor:torus-pairs}. The segment $[-2u,2u]$ corresponds to $r_u=4$, with $D=4$ and $M=4$, and the parallelogram with vertices $\pm u,\pm v$ corresponds to $r_{v-u}=r_{v+u}=1$, with $D=4$ and $M=\Delta$. It remains to show that $M(r)\geq\min\{4,\Delta\}$ for every admissible family $r$ with $D(r)\geq4$. We distinguish cases according to $\sum_er_e$. If $\sum_er_e\leq2$ and $r$ is supported on a single direction, then $D(r)\leq2$, so this case does not occur.
		
		{\bf Case $\sum_er_e\geq4$.} Every primitive lattice vector has length at least $\sys(T^2)=1$, so $M(r)\geq4$.
		
		{\bf Case $\sum_er_e=3$.} By admissibility, $r$ has three distinct directions, each with multiplicity $1$, representing the three nonzero classes of $\Lambda/2\Lambda$.
		
		{\bf Claim.} The shortest representatives of these three classes are $u$, $v$ and $v-u$.
		
		Indeed, a representative $mu+nv$ has length $\sqrt{(m+nx)^2+n^2A^2}\geq|n|A$. In the class of $u$ ($m$ odd, $n$ even), either $n=0$ and the length is at least $1$, or $|n|\geq2$ and the length is at least $2A>1$. In the class of $v$ ($m$ even, $n$ odd), the minimum length is attained at $m=0$, $n=\pm1$, since $|n|\geq3$ implies $\sqrt{(m+nx)^2+n^2A^2}\geq3A>\sqrt{\frac14+A^2}\geq|v|$. In the class of $u+v$ ($m$ and $n$ odd), the minimum length is attained at $m=-1$, $n=1$, since $|n|\geq3$ implies a length of at least $3A>\sqrt{1+A^2}\geq|v-u|$. This proves the claim. Hence, by the strict triangle inequality for the linearly independent vectors $u$ and $v$,
		\[
		M(r) \geq 1+|v|+|v-u| > |v+u|+|v-u|= \Delta .
		\]
		
		{\bf Case $\sum_er_e=2$.} Let $e$ and $f$ be the two distinct directions. Admissibility gives $e\equiv f\pmod{2\Lambda}$. Then $\det(e,f)/A$ is even and nonzero, so $D=2+|\det(e,f)|/A\geq4$. Write $e=m_1u+n_1v$ and $f=m_2u+n_2v$, so that $n_1\equiv n_2\pmod2$, and compare $M(r)=|e|+|f|$ with $\Delta$.
		\begin{itemize}
			\item[($\alpha$)] If $|n_1|+|n_2|\geq4$, then $M(r)\geq(|n_1|+|n_2|)A\geq4A\geq\Delta$ by \eqref{eq:key}.
			\item[($\beta$)] If $n_1$ and $n_2$ are both odd, then by ($\alpha$) we may assume $|n_1|=|n_2|=1$. After changing signs, $e=v+ku$ and $f=v+lu$ with $k>l$. Since $k\equiv l\pmod2$, the difference $n=k-l$ is even and positive, and Lemma \ref{lem:bal} gives $M(r)\geq\Delta$.
			\item[($\gamma$)] If $n_1$ and $n_2$ are both even, then by ($\alpha$) one of them vanishes and the other has absolute value $2$. Thus, after relabeling, $e=\pm u$ and $|n_2|=2$, so $M(r)\geq1+2A\geq\Delta$ by Lemma \ref{lem:key}.
		\end{itemize}
		This completes the proof.
	\end{proof}
	
	\begin{corollary}\label{cor:unsat}
		If $C_3<4$, then $C_3=\Delta$, and $C_3$ determines $x$.
	\end{corollary}
	
	\begin{proof}
		The first assertion follows from Proposition \ref{prop:C3}. By \eqref{eq.deriv.|v+tu|},
		\[
		\partial_x\Delta=\frac{x+1}{\sqrt{(x+1)^2+A^2}}+\frac{x-1}{\sqrt{(x-1)^2+A^2}}=\psi(x+1)+\psi(x-1).
		\]
		Since $\psi$ is odd and strictly increasing, this derivative is positive if and only if $x+1>1-x$, that is, if and only if $x>0$. Hence $\Delta$ is strictly increasing in $x$ on $[0,\frac12]$. Since $A$ is determined by the widths, $C_3=\Delta$ determines $x$.
	\end{proof}
	
	\subsection{Thin tori and hexagons of height two}
	
	For the rest of this section, assume $\Delta\geq4$. Since $\Delta$ is increasing in both $x$ and $A$, and since $\Delta=(\sqrt{109}+\sqrt{181})/6<4$ for $x=\frac12$ and $A=\frac53$, we must have $A>\frac53$, so the torus is elongated.
	
	Write each direction as $e=m_eu+n_ev$ with $m_e,n_e\in\Z$, and define the {\bf height} of an admissible family $r$ (or of the polygon $Z(r)$) to be $\sum_er_e|n_e|$. The vertical width of $Z(r)$ equals $A$ times its height. Since the perimeter of a convex set is at least twice its width in any direction,
	\begin{equation}\label{eq:hbound}
		M(r)\geq A \cdot \text{height}.
	\end{equation}
	By admissibility, $\sum_er_en_e$ is even, and hence so is the height. Consequently, every admissible family with $M(r)<4A$ has height $0$ or $2$. Polygons of height $0$ are horizontal segments.
	
	A polygon of height $2$ has horizontal multiplicity $\rho\geq0$ and either two nonhorizontal generators $v+ku$, $v+lu$ with $k\geq l$, counted with multiplicity, or one generator $ku+2v$. In the first case, write $n=k-l$. The edge vectors sum to $(\rho+k+l)u+2v$, so admissibility gives $\rho\equiv n\pmod2$, and
	\begin{equation}\label{eq:hexdata}
		D=\underbrace{(\rho+2)}_{\sum_e r_e}+\underbrace{2\rho+n}_{\sum_{e<f} r_er_f|\det(e,f)|/A}
		=3\rho+n+2,\qquad M=\rho+|v+ku|+|v+lu| .
	\end{equation}
	We note two consequences of \eqref{eq:hexdata}. First, increasing $\rho$ by $2$ changes $(D,M)$ by $(6,2)$, and the quantity
	\begin{equation}\label{eq:betakl}
		M-D/3=|v+ku|+|v+lu|-\tfrac{n+2}{3}
	\end{equation}
	is independent of $\rho$. Second, $\rho\equiv n\pmod2$ gives $3\rho\equiv3n\pmod6$, and therefore
	\begin{equation}\label{eq.D.mod6}
		D \equiv 4n+2 \pmod 6.
	\end{equation}
	Thus the hexagons split into three families, motivating the following definition.
	
	\begin{definition}\label{def:beta}
		For integers $k\geq l$, define
		\begin{equation}\label{eq:beta}
			\beta(k,l):=|v+ku|+|v+lu|-\frac{k-l+2}{3},
		\end{equation}
		and for $j\in\Z/3$, define
		\begin{equation}\label{eq:betaj}
			\beta_j:=\min\bigl\{\beta(k,l)\ :\ k\geq l,\ k-l\equiv j\!\!\pmod 3\bigr\} .
		\end{equation}
	\end{definition}
	
	By \eqref{eq:betakl} and \eqref{eq.D.mod6}, $\beta_j$ is the least value of $M-\frac13D$ among the hexagons whose lattice count $D$ is congruent to $4j+2$ modulo $6$. Equivalently, $M=\frac13D+\beta_j$ is the lowest line of slope $\frac13$ in the $(D,M)$-plane that meets this family.
	
	\begin{lemma}\label{lem:beta}
		The numbers $\beta_0,\beta_1,\beta_2$ have the following properties.
		\begin{enumerate}
			\item[(i)] $\beta(k,l)\geq\frac{2(k-l)-2}{3}$. In particular, the minimum in \eqref{eq:betaj} is attained, and every minimizer satisfies $k-l\leq\frac32\beta_j+1$.
			\item[(ii)] $\beta_j<2A-\frac12$.
			\item[(iii)] $\beta_{j+1}\leq\beta_j+\frac23$, and hence $\beta_{j+r}\leq\beta_j+\frac{2r}{3}$ for every integer $r\geq0$.
		\end{enumerate}
	\end{lemma}
	
	\begin{proof}
		(i) Projecting to the horizontal axis gives $|v+ku|+|v+lu|\geq|x+k|+|x+l|\geq k-l$, so $\beta(k,l)\geq(k-l)-\frac{(k-l)+2}{3}=\frac{2(k-l)-2}{3}$. Since $\beta_j\leq\beta(j,0)<\infty$, the relevant competitors in \eqref{eq:betaj} have bounded $k-l$, and for each value of $k-l$, only finitely many pairs satisfy $\beta(k,l)\leq\beta(j,0)$. Hence the minimum is attained, and every minimizer satisfies $\frac{2(k-l)-2}{3}\leq\beta_j$.
		
		(ii) The pairs $(k,l)=(0,0),\,(0,-1),\,(1,-1)$, with differences $0,1,2$, give $\beta_0\leq2|v|-\frac23$, $\beta_1\leq|v|+|v-u|-1$ and $\beta_2\leq\Delta-\frac43$. Using $\sqrt{(x+t)^2+A^2}\leq A+\frac{(x+t)^2}{2A}$ and $0\leq x\leq\frac12$, we obtain
		\[
		2|v|\leq2A+\frac{1}{4A},\qquad
		|v|+|v-u|\leq2A+\frac{1}{2A},\qquad
		\Delta\leq2A+\frac{5}{4A}.
		\]
		Since $A>\frac53$, each of $2A+\frac{1}{4A}-\frac23$, $2A+\frac{1}{2A}-1$ and $2A+\frac{5}{4A}-\frac43$ is smaller than $2A-\frac12$.
		
		(iii) Suppose that $\beta_j$ is attained at $(k,l)$. Then $(k+1)-l\equiv j+1\pmod3$, and
		\[
		\beta(k+1,l)-\beta(k,l)=|v+(k+1)u|-|v+ku|-\frac13\leq1-\frac13=\frac23.
		\]
		The second assertion follows by induction.
	\end{proof}
	
	We now bound $M$ according to the height of the polygon.
	
	\begin{proposition}\label{prop:classify}
		Let $r$ be a nonzero admissible family. Then its height is even, and the following hold.
		\begin{enumerate}
			\item[(a)] If the height is $0$, then $Z(r)$ is a horizontal segment $[-iu,iu]$, and $D(r)=M(r)=2i$.
			\item[(b)] If the height is $2$, then $D(r)\equiv4j+2\pmod6$ for a unique $j\in\Z/3$, and
			\begin{equation}\label{eq:aboveline}
				M(r)\geq \tfrac13D(r)+\beta_j .
			\end{equation}
			\item[(c)] If the height is at least $4$, then $M(r)\geq 4A$.
		\end{enumerate}
		Moreover, for each $j\in\Z/3$, there exists $d_j\equiv4j+2\pmod6$ with $d_j\leq\frac32\beta_j+6$ such that
		\begin{equation}\label{eq:attain}
			C_{d_j+6t}\ \leq\ \tfrac13(d_j+6t)+\beta_j\qquad\text{for every integer }t\geq0 .
		\end{equation}
	\end{proposition}
	
	\begin{proof}
		Part (c) follows from \eqref{eq:hbound}.
		
		(a) Every direction in the support of $r$ has $n_e=0$, and hence equals $u$. By admissibility, $r_u=2i$ is even, so $D(r)=M(r)=2i$.
		
		(b) Suppose first that the nonhorizontal generators are $v+ku$ and $v+lu$. Then \eqref{eq:hexdata} holds, and \eqref{eq.D.mod6} identifies $j$ with $n\bmod3$. Hence, by \eqref{eq:betakl} and \eqref{eq:betaj},
		\[
		M-\frac13 D=\beta(k,l)\geq \beta_j.
		\]
		If the only nonhorizontal generator is $ku+2v$, then $k$ is odd. Admissibility then implies that $\rho$ is odd, and since $|\det(u,ku+2v)|/A=2$, we obtain $D=(\rho+1)+2\rho=3\rho+1\equiv4\pmod6$, so $j=2$. Moreover, $M=\rho+|ku+2v|\geq\rho+|2v-u|$, because the minimum over odd $k$ is attained at $k=-1$ as $0\leq 2x\leq 1$. Hence
		\[
		M-\frac13D\geq |2v-u|-\frac13\geq \Delta-\frac43\geq \beta_2,
		\]
		where the second inequality follows from \eqref{eq:key}, and the last one from the choice $(k,l)=(1,-1)$ in \eqref{eq:betaj}.
		
		For the last assertion, let $(k,l)$ attain $\beta_j$ and let $n=k-l$. Consider the hexagon in \eqref{eq:hexdata} with these generators and with $\rho$ the least nonnegative integer of the same parity as $n$. Increasing $\rho$ by $2$ repeatedly, we realize every lattice count in $d_j+6\Z_{\geq0}$ with $M=\frac13D+\beta_j$, where $d_j=n+2$ if $n$ is even and $d_j=n+5$ if $n$ is odd. In either case, Lemma \ref{lem:beta}(i) gives $d_j\leq n+5\leq\frac32\beta_j+6$.
	\end{proof}
	
	The next proposition recovers each $\beta_j$ from the sequence of widths.
	
	\begin{proposition}\label{prop.beta.C_p}
		For every $j\in\Z/3$,
		\begin{equation}\label{eq:beta.C_p}
			\beta_j=\max\Bigl\{C_p-\frac p3:\ p\equiv 4j+2\pmod{6},\ C_p<4A\Bigr\}.
		\end{equation}
	\end{proposition}
	
	\begin{proof}
		Fix $j$ and write $\beta=\beta_j$. Choose a pair $(k,l)$ attaining $\beta$, let $n=k-l$, and choose $d=d_j$ as in Proposition \ref{prop:classify}, so that $d\leq n+5$ and $d\leq\frac32\beta+6$.
		
		We first prove that the right-hand side of \eqref{eq:beta.C_p} is at most $\beta$. Let $p\equiv4j+2\pmod{6}$ with $C_p<4A$; in particular, $p$ is even. If $p\geq d$, then $p\in d+6\Z_{\geq0}$, and \eqref{eq:attain} gives $C_p\leq\frac p3+\beta$. If $p<d$, then $p\leq d-6\leq n-1$. The segment $[-\frac p2u,\frac p2u]$ gives $C_p\leq p$, so $C_p-\frac p3\leq\frac{2p}{3}\leq\frac{2n-2}{3}\leq\beta$ by Lemma \ref{lem:beta}(i).
		
		We next prove the reverse inequality. It suffices to find $p\in d+6\Z_{\geq0}$ in the range
		\begin{equation}\label{eq.range.p}
			\tfrac32\beta\ \leq\ p\ <\ 12A-3\beta .
		\end{equation}
		Indeed, suppose that $p$ lies in this range, and let $r$ be an admissible family with $D:=D(r)\geq p$. We verify $M(r)\geq\frac p3+\beta$ in each of the three cases of Proposition \ref{prop:classify}.
		
		If the height is $0$, then $M=D\geq p\geq\frac p3+\beta$, where the last inequality follows from the lower bound in \eqref{eq.range.p}.
		
		If the height is at least $4$, then $M\geq4A>\frac p3+\beta$ by the upper bound in \eqref{eq.range.p}.
		
		If the height is $2$, write $D=p+2c+6t$ with $c\in\{0,1,2\}$ and an integer $t\geq0$. Comparing \eqref{eq.D.mod6} for $D$ with $p\equiv4j+2\pmod6$ gives $4n'\equiv4j+2c\pmod6$ for the corresponding $n'$, that is, $n'\equiv j+2c\pmod3$. Hence \eqref{eq:aboveline} and Lemma \ref{lem:beta}(iii) give
		\[
		M\geq \frac D3+\beta_{j+2c}
		=\frac p3+\frac{2c}{3}+2t+\beta_{j+2c}\geq \frac p3+\beta,
		\]
		where the last inequality follows from $\beta=\beta_{(j+2c)+c}\leq\beta_{j+2c}+\frac{2c}{3}$.
		
		Therefore $C_p\geq\frac p3+\beta$, and \eqref{eq:attain} gives the reverse inequality. Thus $C_p=\frac p3+\beta<4A$ by the upper bound in \eqref{eq.range.p}, so $p$ belongs to the set in \eqref{eq:beta.C_p}, and the maximum is at least $\beta$.
		
		Finally, we show that such a $p$ exists. By Lemma \ref{lem:beta}(ii) and $A>\frac53$, the interval \eqref{eq.range.p} has length $12A-\frac92\beta>3A+\frac94>6$. The progression $d+6\Z_{\geq0}$ has spacing $6$ and begins at $d\leq\frac32\beta+6$. Hence the first of them that is at least $\frac32\beta$ is at most $\frac32\beta+6$, and therefore lies in the interval \eqref{eq.range.p}.
	\end{proof}
	
	\subsection{Recovery of the shear parameter}
	
	We show that $\beta_0,\beta_1,\beta_2$ determine $x$.
	
	\begin{lemma}\label{lem:monotone}
		On $[0,\frac12]$, the function $x\to\beta(k,l)$ is strictly increasing if $k+l\geq0$, and nonincreasing if $k+l\leq-1$. Moreover, at $x=\frac12$, at least one of $\beta_0,\beta_1,\beta_2$ is attained at a pair $(k,l)$ with $k+l\geq0$.
	\end{lemma}
	
	\begin{proof}
		\textbf{Monotonicity.} By \eqref{eq.deriv.|v+tu|}, $\partial_x\beta(k,l)=\psi(x+k)+\psi(x+l)$. Since $\psi$ is odd and strictly increasing, $\partial_x\beta(k,l)>0$ if and only if $x+k>-(x+l)$, that is, if and only if $2x+k+l>0$. If $k+l\geq0$, this holds on $(0,\frac12]$; if $k+l\leq -1$, it fails on $[0,\frac12]$.
		
		\textbf{Attainment.} By Lemma \ref{lem:bal}, for fixed $n=k-l$, the minimum of $\beta(k,l)$ is attained at the balanced pair, namely $(m,-m)$ if $n=2m$ and $(m,-m-1)$ if $n=2m+1$. These pairs satisfy $k+l=0$ and $k+l=-1$, respectively. It therefore suffices to show that the three minima cannot all occur only at odd values of $n$.
		
		Let $\lambda_m=\sqrt{(m+\frac12)^2+A^2}$, so that at $x=\frac12$ we have $|v+mu|=\lambda_m$ and $|v-mu|=\lambda_{m-1}$, with $\lambda_{-1}=\lambda_0$. Denote by $b_n$ the value at $x=\frac12$ of $\beta$ at the balanced pair with $k-l=n$:
		\[
		b_{2m}=\lambda_m+\lambda_{m-1}-\tfrac{2m+2}{3},\qquad
		b_{2m+1}=2\lambda_m-\tfrac{2m+3}{3}.
		\]
		With $a_m:=\lambda_m-\lambda_{m-1}-\frac13$, we have
		\begin{equation}\label{eq:increments}
			b_{2m+1}-b_{2m}=a_m,\qquad b_{2m+2}-b_{2m+1}=a_{m+1} ,
		\end{equation}
		and the strict convexity of $\tau\to\sqrt{\tau^2+A^2}$ implies that $(a_m)_{m\geq0}$ is strictly increasing. Moreover, $\beta_j(\frac12)=\min\{b_n:n\equiv j \bmod 3\}$.
		
		Suppose that all three minima are attained only at odd indices $n_i=2m_i+1$, $i=1,2,3$, with $m_1<m_2<m_3$. Since the $n_i$ lie in distinct classes modulo $3$ and $m\to2m+1$ is a bijection of $\Z/3$, the $m_i$ are also distinct modulo $3$, and hence $m_3\geq m_1+2$. By \eqref{eq:increments},
		\[
		b_{2m+4}-b_{2m+1}=2a_{m+1}+a_{m+2},\qquad
		b_{2m+1}-b_{2m-2}=a_{m-1}+2a_m,
		\]
		and both indices $2m+4$ and $2m-2$ lie in the same class modulo $3$ as $2m+1$. The minimality at $n_1$ and at $n_3$ gives
		\[
		2a_{m_1+1}+a_{m_1+2}\geq0,\qquad a_{m_3-1}+2a_{m_3}\leq0 .
		\]
		Since $m_1+1\leq m_3-1$, $m_1+2\leq m_3$, and $(a_m)$ is increasing, the second inequality yields $a_{m_1+1}+2a_{m_1+2}\leq0$. Subtracting this from the first inequality gives $a_{m_1+1}-a_{m_1+2}\geq0$, which contradicts the strict monotonicity of $(a_m)$.
	\end{proof}
	
	The preceding lemma yields the required monotonicity.
	
	\begin{proposition}\label{prop:xdet}
		For each $A>\frac53$, there exists $j^*\in\Z/3$, depending only on $A$, such that the function $x\to\beta_{j^*}$ is strictly increasing on $[0,\frac12]$. In particular, the triple $(\beta_0,\beta_1,\beta_2)$ determines $x$.
	\end{proposition}
	
	\begin{proof}
		Throughout the proof, $\beta(k,l)$ and $\beta_j$ are regarded as functions of $x\in[0,\frac12]$. By Lemma \ref{lem:monotone}, we may choose $j^*$ and a pair $(k^*,l^*)$ with $k^*-l^*\equiv j^*\pmod3$ and $k^*+l^*\geq0$ that attains $\beta_{j^*}$ at $x=\frac12$. If $(k,l)$ is any competitor in the same class with $k+l\leq -1$, then for every $x\in[0,\frac12]$,
		\[
		\beta(k^*,l^*)\big|_x\leq \beta(k^*,l^*)\big|_{1/2}\leq
		\beta(k,l)\big|_{1/2}\leq \beta(k,l)\big|_x ,
		\]
		where the first and third inequalities follow from Lemma \ref{lem:monotone}, and the second follows from the choice of $(k^*,l^*)$. Hence, for every $x$, the minimum defining $\beta_{j^*}$ is attained at some pair with $k+l\geq0$. Since the minimum is attained, it is strictly increasing: if $x_1<x_2$ and $(k_2,l_2)$ attains $\beta_{j^*}$ at $x_2$ with $k_2+l_2\geq 0$, then $\beta_{j^*}\big|_{x_2}=\beta(k_2,l_2)\big|_{x_2} >\beta(k_2,l_2)\big|_{x_1}\geq\beta_{j^*}\big|_{x_1}$.
	\end{proof}
	
	\begin{proof}[Proof of Theorem \ref{thm:flat_rigidity}]
		Let $T$ and $T'$ be flat two-dimensional tori with $\omega_p(T)=\omega_p(T')$ for all $p$. Lemma \ref{lem:1stw_sys} gives $\sys(T)=\sys(T')$, and \eqref{eq:Weyl} gives $\area(T)=\area(T')$. Rescale both tori to have systole $1$. Their widths still agree, and their lattices have the same area $A$. If $C_3<4$, then Corollary \ref{cor:unsat} gives the same shear parameter $x$ for both lattices. Otherwise, $A>\frac53$, Proposition \ref{prop.beta.C_p} recovers the same triple $(\beta_0,\beta_1,\beta_2)$ from the common widths, and Proposition \ref{prop:xdet} shows that the two lattices have the same value of $x$. Thus the lattices, and hence the tori, are isometric.
	\end{proof}
	
	\begin{remark}
		Both cases occur.
		\begin{itemize}
			\item For the square lattice, $x=0$ and $A=1$, we have $\Delta=2\sqrt2<4$, so the shear parameter is already detected by the third width.
			\item For the lattice with $x=0$ and $A=2$, we have $\Delta=2\sqrt5>4$, and the hexagons of height two are needed.
		\end{itemize}
	\end{remark}

	\section{Classification of the quotient of flat tori}\label{S:quotient classification}
	
	Theorem \ref{thm:affine} applies to every isometric quotient of flat tori, which is also a compact flat $2$-orbifold. In this appendix we recall the classification of Conway and Huson \cite{ConwayHuson02}, and illustrate the seventeen types.

	A compact flat $2$-orbifold is a quotient of $\R^2$ by a discrete cocompact group $\Gamma$ of Euclidean isometries. Its translation subgroup is a lattice of finite index, so every such orbifold can be written as a finite isometric quotient $T^2/G$ of a flat torus. Conversely, every finite isometric quotient of a flat torus is a compact flat $2$-orbifold. Up to affine conjugacy, there are seventeen types. We use the orbifold notation of Conway and Huson \cite{ConwayHuson02}: $\circ$ denotes a handle, $\times$ a cross-cap, and $\ast$ a boundary component; integers before $\ast$ record cone orders and integers after $\ast$ record corner orders.
	
	\begin{table}[tbp]
		\centering
		\small
		\begin{tabular}{@{}llll@{}}
			\toprule
			Signature & Underlying surface & Singular data & Crystallographic symbol\\
			\midrule
			$\circ$ & torus & none & \texttt{p1}\\
			$2222$ & sphere & four cone points of order $2$ & \texttt{p2}\\
			$333$ & sphere & three cone points of order $3$ & \texttt{p3}\\
			$442$ & sphere & cone orders $4,4,2$ & \texttt{p4}\\
			$632$ & sphere & cone orders $6,3,2$ & \texttt{p6}\\
			\midrule
			$\ast\ast$ & annulus & mirror boundary & \texttt{pm}\\
			$\ast\times$ & M\"obius band & mirror boundary & \texttt{cm}\\
			$\times\times$ & Klein bottle & none & \texttt{pg}\\
			\midrule
			$\ast2222$ & disk & four corners of order $2$ & \texttt{pmm}\\
			$2{\ast}22$ & disk & one cone point of order $2$, two corners of order $2$ & \texttt{cmm}\\
			$22{\ast}$ & disk & two cone points of order $2$ & \texttt{pmg}\\
			$22\times$ & projective plane & two cone points of order $2$ & \texttt{pgg}\\
			\midrule
			$\ast333$ & disk & three corners of order $3$ & \texttt{p3m1}\\
			$3{\ast}3$ & disk & one cone point and one corner of order $3$ & \texttt{p31m}\\
			$\ast442$ & disk & corner orders $4,4,2$ & \texttt{p4m}\\
			$4{\ast}2$ & disk & one cone point of order $4$, one corner of order $2$ & \texttt{p4g}\\
			$\ast632$ & disk & corner orders $6,3,2$ & \texttt{p6m}\\
			\bottomrule
		\end{tabular}
		\caption{The seventeen compact flat $2$-orbifolds in Conway--Huson notation \cite{ConwayHuson02}. A cone point of order $m$ has angle $2\pi/m$, while a corner of order $m$ has angle $\pi/m$.}
		\label{quot:tab:types}
	\end{table}
	
	The following figures give schematic fundamental domains for all seventeen types. Dashed edges are identified according to the arrows, thick edges are mirror boundary, circles mark cone points, and squares mark mirror corners. The pictures indicate the topology and singular data; the Euclidean moduli are not encoded by the drawings.
	
	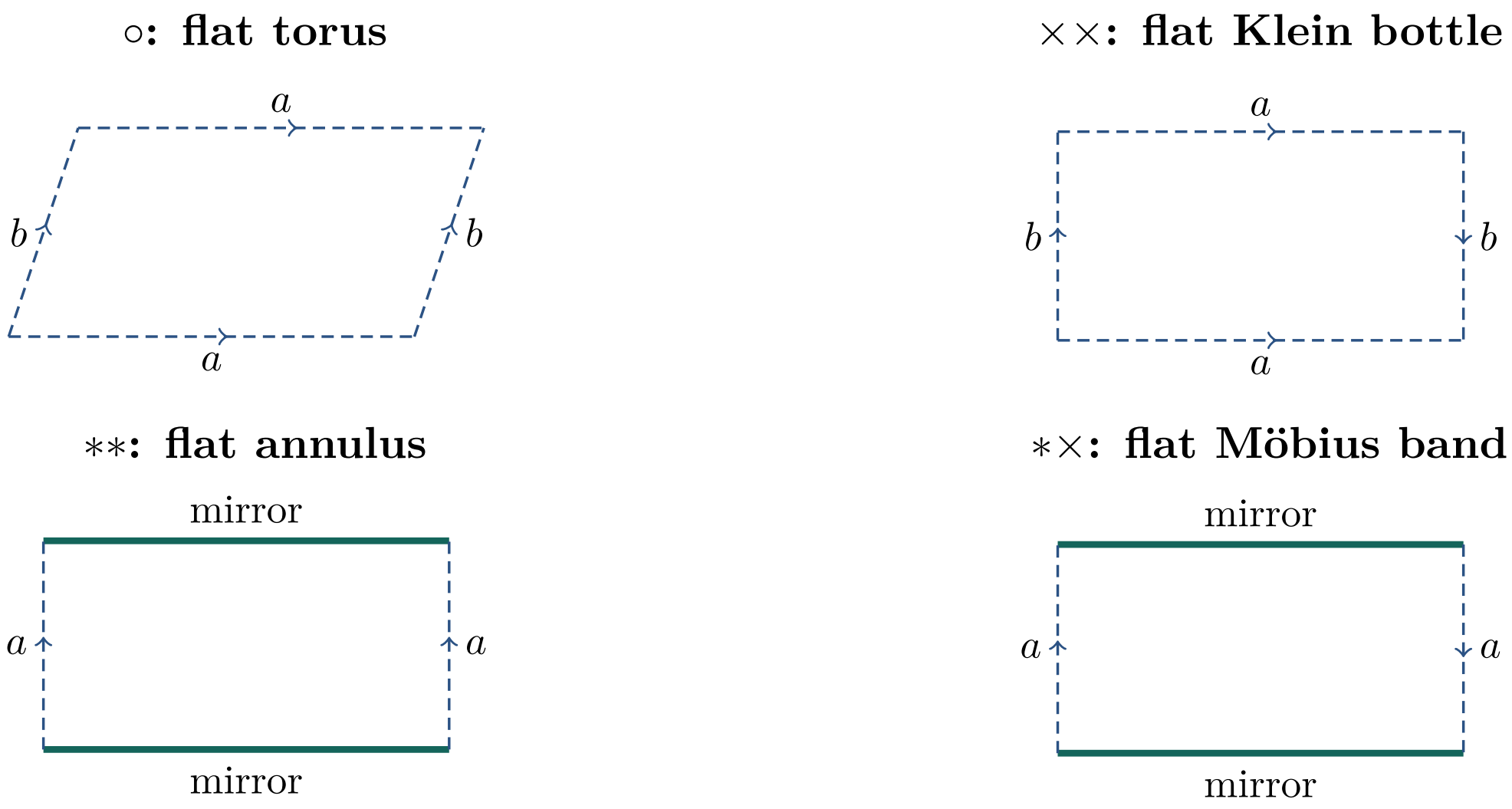
\begin{figure}[tbp]
		\centering
		\newcommand{\qrect}[5]{
			\begin{minipage}[t]{.47\textwidth}\centering\textbf{#1}\par\smallskip
				\begin{tikzpicture}[font=\small,x=1cm,y=1cm]
					\path[use as bounding box] (-.55,-.55) rectangle (4.8,2.4);#2#3#4#5
		\end{tikzpicture}\end{minipage}}
		\qrect{$\circ$: flat torus}
		{\draw[qcut] (0,0) -- node[below] {$a$} (3.5,0);}
		{\draw[qcut] (.6,1.8) -- node[above] {$a$} (4.1,1.8);}
		{\draw[qcut] (0,0) -- node[left] {$b$} (.6,1.8);}
		{\draw[qcut] (3.5,0) -- node[right] {$b$} (4.1,1.8);}
		\hfill
		\qrect{$\times\times$: flat Klein bottle}
		{\draw[qcut] (.3,0) -- node[below] {$a$} (3.8,0);}
		{\draw[qcut] (.3,1.8) -- node[above] {$a$} (3.8,1.8);}
		{\draw[qcut] (.3,0) -- node[left] {$b$} (.3,1.8);}
		{\draw[qcut] (3.8,1.8) -- node[right] {$b$} (3.8,0);}
		\par\medskip
		\qrect{$\ast\ast$: flat annulus}
		{\draw[qmirror] (.3,0) -- node[below] {mirror} (3.8,0);}
		{\draw[qmirror] (.3,1.8) -- node[above] {mirror} (3.8,1.8);}
		{\draw[qcut] (.3,0) -- node[left] {$a$} (.3,1.8);}
		{\draw[qcut] (3.8,0) -- node[right] {$a$} (3.8,1.8);}
		\hfill
		\qrect{$\ast\times$: flat M\"obius band}
		{\draw[qmirror] (.3,0) -- node[below] {mirror} (3.8,0);}
		{\draw[qmirror] (.3,1.8) -- node[above] {mirror} (3.8,1.8);}
		{\draw[qcut] (.3,0) -- node[left] {$a$} (.3,1.8);}
		{\draw[qcut] (3.8,1.8) -- node[right] {$a$} (3.8,0);}
		\caption{The four compact flat orbifolds without cone points or corners.}
		\label{fig:quotient-smooth}
	\end{figure}
	
	\begin{figure}[tbp]
		\centering
		\newcommand{\qdouble}[2]{
			\begin{tikzpicture}[font=\small,x=1cm,y=1cm]
				\path[use as bounding box] (-.65,-.55) rectangle (5.95,2.15);
				\foreach \sh in {0,3.35}{\begin{scope}[shift={(\sh,0)}]#1#2\end{scope}}
		\end{tikzpicture}}
		\begin{minipage}[t]{.47\textwidth}\centering
			\textbf{$2222$}\par\smallskip
			\qdouble{\draw[qcut] (0,0) -- node[below] {$a$} (1.75,0);
				\draw[qcut] (1.75,0) -- node[right] {$b$} (1.75,1.2);
				\draw[qcut] (1.75,1.2) -- node[above] {$c$} (0,1.2);
				\draw[qcut] (0,1.2) -- node[left] {$d$} (0,0);}
			{\qpt{qcone}{(0,0)}{below left}{2}\qpt{qcone}{(1.75,0)}{below right}{2}
				\qpt{qcone}{(1.75,1.2)}{above right}{2}\qpt{qcone}{(0,1.2)}{above left}{2}}
		\end{minipage}\hfill
		\begin{minipage}[t]{.47\textwidth}\centering
			\textbf{$333$}\par\smallskip
			\qdouble{\draw[qcut] (0,0) -- node[below] {$a$} (1.8,0);
				\draw[qcut] (1.8,0) -- node[right] {$b$} (0.9,1.5588);
				\draw[qcut] (0.9,1.5588) -- node[left] {$c$} (0,0);}
			{\qpt{qcone}{(0,0)}{below left}{3}\qpt{qcone}{(1.8,0)}{below right}{3}\qpt{qcone}{(0.9,1.5588)}{above}{3}}
		\end{minipage}
		\par\medskip
		\begin{minipage}[t]{.47\textwidth}\centering
			\textbf{$442$}\par\smallskip
			\qdouble{\draw[qcut] (0,0) -- node[below] {$a$} (1.65,0);
				\draw[qcut] (1.65,0) -- node[above right] {$b$} (0,1.65);
				\draw[qcut] (0,1.65) -- node[left] {$c$} (0,0);}
			{\qpt{qcone}{(0,0)}{below left}{2}\qpt{qcone}{(1.65,0)}{below right}{4}\qpt{qcone}{(0,1.65)}{above}{4}}
		\end{minipage}\hfill
		\begin{minipage}[t]{.47\textwidth}\centering
			\textbf{$632$}\par\smallskip
			\qdouble{\draw[qcut] (0,0) -- node[below] {$a$} (2,0);
				\draw[qcut] (2,0) -- node[above right] {$b$} (0,1.1547);
				\draw[qcut] (0,1.1547) -- node[left] {$c$} (0,0);}
			{\qpt{qcone}{(0,0)}{below left}{2}\qpt{qcone}{(2,0)}{below right}{6}\qpt{qcone}{(0,1.1547)}{above}{3}}
		\end{minipage}
		\caption{The four cone spheres. Each is obtained by gluing two congruent polygons along corresponding dashed edges.}
		\label{fig:quotient-cones}
	\end{figure}
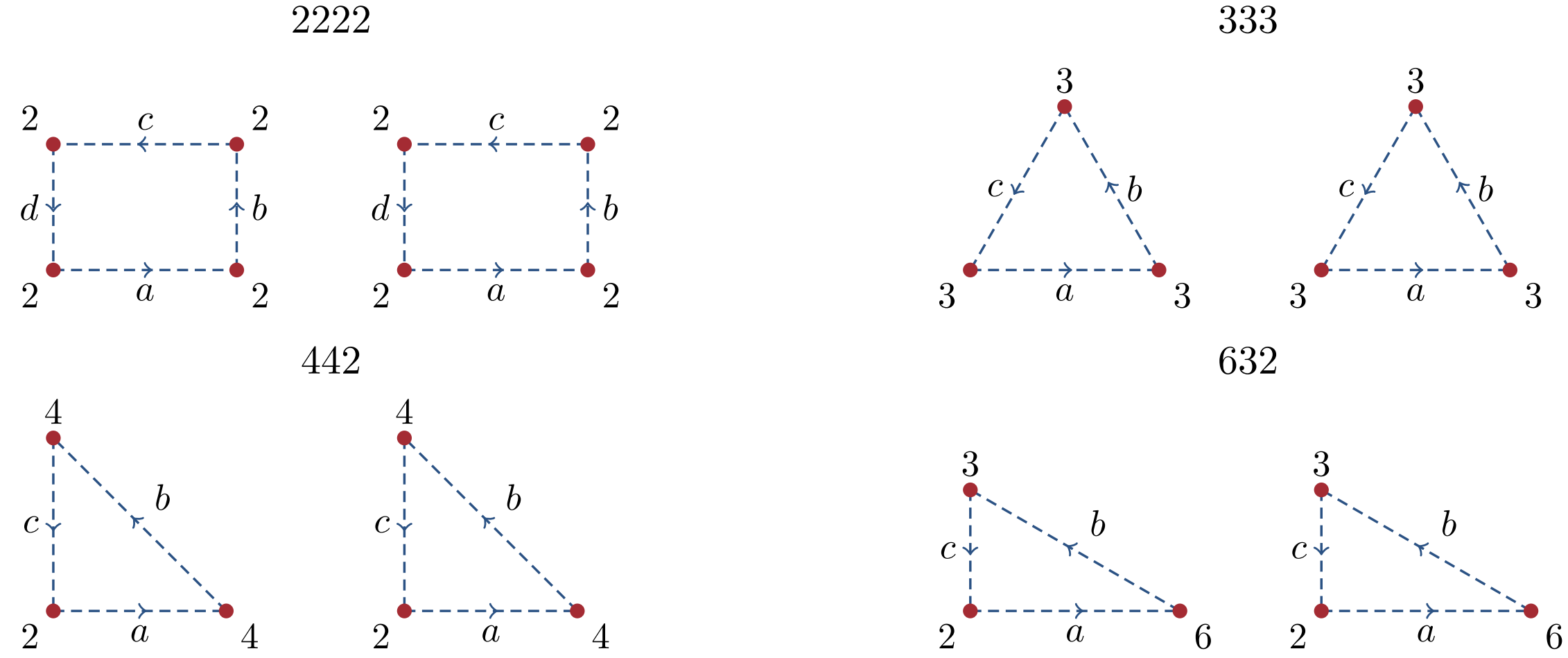
	
	\begin{figure}[tbp]
		\centering
		\newcommand{\qbox}[2]{\begin{minipage}[t]{.47\textwidth}\centering\textbf{#1}\par\smallskip
				\begin{tikzpicture}[font=\small,x=1cm,y=1cm]
					\path[use as bounding box] (-.45,-.45) rectangle (4.55,3.1);#2\end{tikzpicture}\end{minipage}}
		\qbox{$\ast2222$: rectangle}{\draw[qmirror] (.5,.35) rectangle (3.6,2.1);
			\qpt{qcorner}{(.5,.35)}{below left}{2}\qpt{qcorner}{(3.6,.35)}{below right}{2}
			\qpt{qcorner}{(3.6,2.1)}{above right}{2}\qpt{qcorner}{(.5,2.1)}{above left}{2}}
		\hfill
		\qbox{$\ast333$: equilateral triangle}{\draw[qmirror] (.5,.15) -- (3.5,.15) -- (2,2.748) -- cycle;
			\qpt{qcorner}{(.5,.15)}{below left}{3}\qpt{qcorner}{(3.5,.15)}{below right}{3}\qpt{qcorner}{(2,2.748)}{above}{3}}
		\par\medskip
		\qbox{$\ast442$: right isosceles triangle}{\draw[qmirror] (.65,.15) -- (3.15,.15) -- (.65,2.65) -- cycle;
			\qpt{qcorner}{(.65,.15)}{below left}{2}\qpt{qcorner}{(3.15,.15)}{below right}{4}\qpt{qcorner}{(.65,2.65)}{above}{4}}
		\hfill
		\qbox{$\ast632$: $30$--$60$--$90$ triangle}{\draw[qmirror] (.3,.15) -- (4.1,.15) -- (.3,2.344) -- cycle;
			\qpt{qcorner}{(.3,.15)}{below left}{2}\qpt{qcorner}{(4.1,.15)}{below right}{6}\qpt{qcorner}{(.3,2.344)}{above}{3}}
		\caption{The four reflection polygons. The entire outline is mirror boundary.}
		\label{fig:quotient-polygons}
	\end{figure}
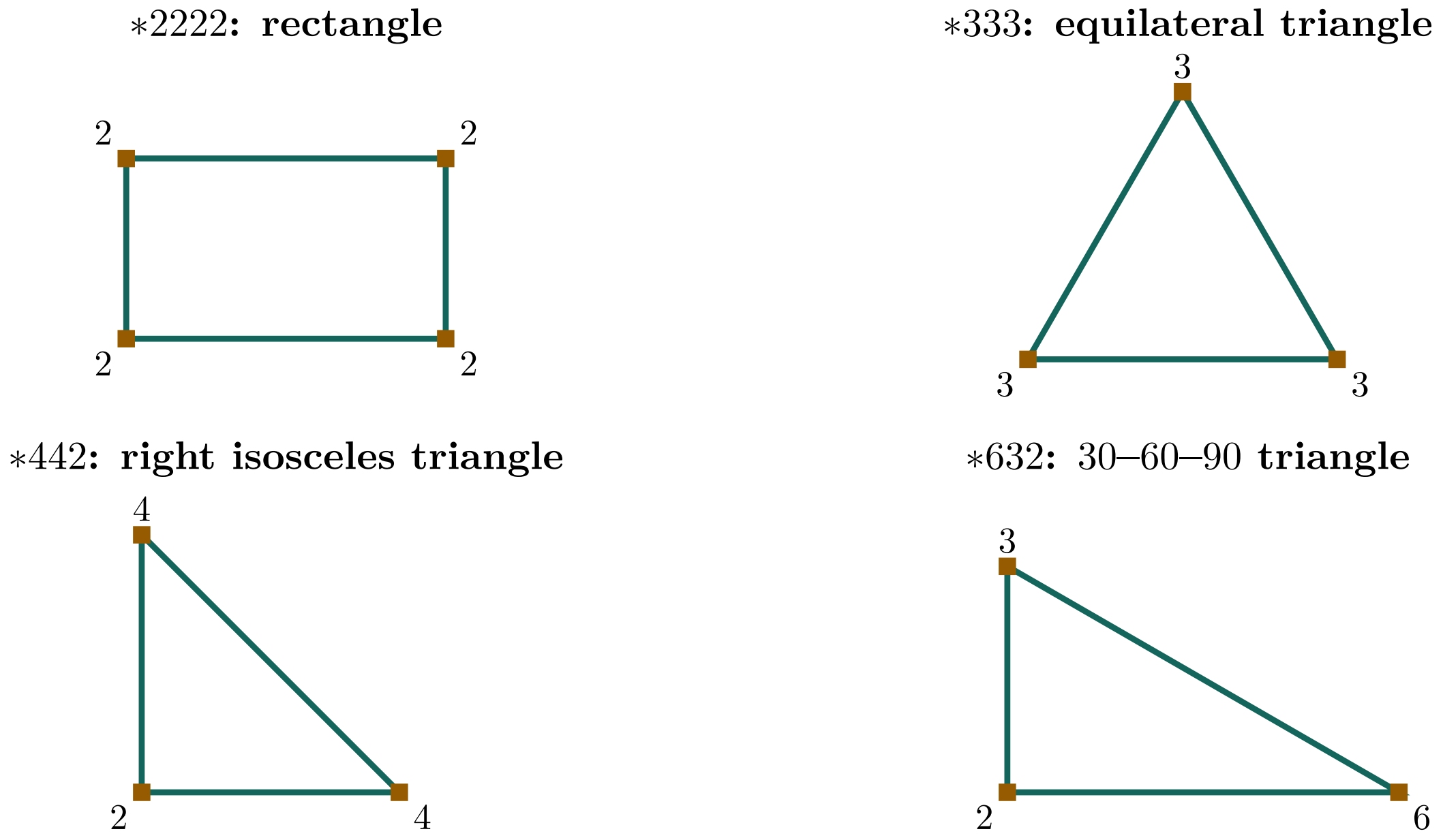
	
	\begin{figure}[tbp]
		\centering
		\newcommand{\qdisk}[2]{\begin{minipage}[t]{.47\textwidth}\centering\textbf{#1}\par\smallskip
				\begin{tikzpicture}[font=\small,x=1cm,y=1cm]
					\path[use as bounding box] (-2.25,-1.45) rectangle (2.25,1.5);#2\end{tikzpicture}\end{minipage}}
		\qdisk{$2{\ast}22$}{\draw[qmirror] (-1.65,0) .. controls (-.7,.95) and (.7,.95) .. (1.65,0)
			.. controls (.7,-.95) and (-.7,-.95) .. (-1.65,0);
			\qpt{qcorner}{(-1.65,0)}{left}{2}\qpt{qcorner}{(1.65,0)}{right}{2}\qpt{qcone}{(0,0)}{above}{2}}
		\hfill
		\qdisk{$3{\ast}3$}{\draw[qmirror] (1.65,0) .. controls (.7,.55) and (-1.65,1.4) .. (-1.65,0)
			.. controls (-1.65,-1.4) and (.7,-.55) .. (1.65,0);
			\qpt{qcorner}{(1.65,0)}{right}{3}\qpt{qcone}{(-.4,0)}{above}{3}}
		\par\medskip
		\qdisk{$4{\ast}2$}{\draw[qmirror] (1.65,0) .. controls (.7,.95) and (-1.65,1.4) .. (-1.65,0)
			.. controls (-1.65,-1.4) and (.7,-.95) .. (1.65,0);
			\qpt{qcorner}{(1.65,0)}{right}{2}\qpt{qcone}{(-.4,0)}{above}{4}}
		\hfill
		\qdisk{$22{\ast}$}{\draw[qmirror] (0,0) ellipse [x radius=1.65,y radius=1.02];
			\qpt{qcone}{(-.65,0)}{above}{2}\qpt{qcone}{(.65,0)}{above}{2}}
		\caption{The four disks with interior cone points and, possibly, mirror corners. The curved outlines indicate only the topology and singular data.}
		\label{fig:quotient-mixed}
	\end{figure}
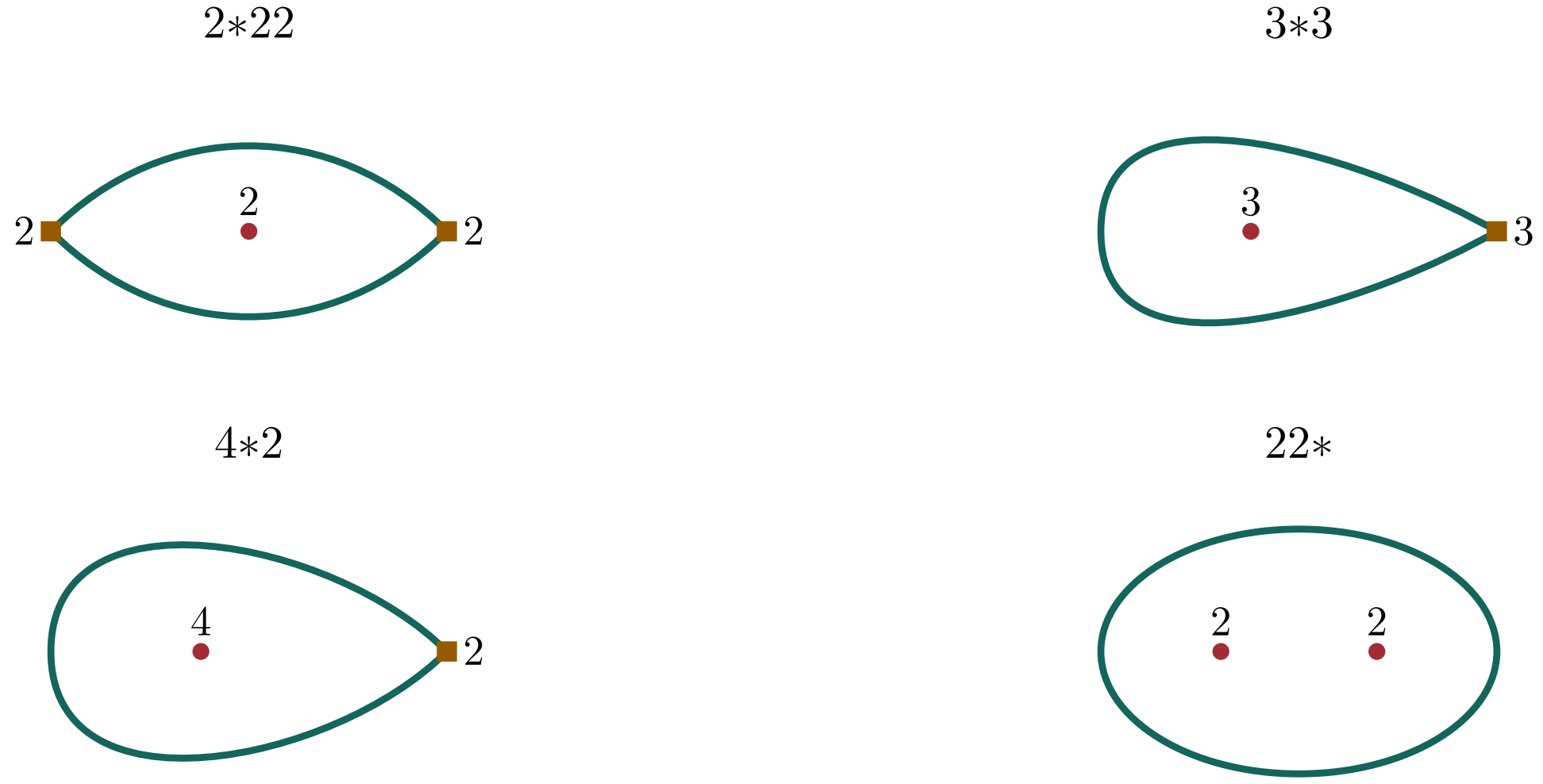
	
	\begin{figure}[tbp]
		\centering
		\begin{tikzpicture}[font=\small,x=1cm,y=1cm]
			\draw[qcut] (1.65,0) arc[start angle=0,end angle=180,radius=1.65];
			\draw[qcut] (-1.65,0) arc[start angle=180,end angle=360,radius=1.65];
			\node[above right] at (1.17,1.17) {$a$};\node[below left] at (-1.17,-1.17) {$a$};
			\foreach \p/\l/\q in {(1.65,0)/right/P,(-1.65,0)/left/P,(0,1.65)/above/Q,(0,-1.65)/below/Q}
			{\draw[quotientglue,fill=white,line width=.65pt] \p circle[radius=1.7pt];\node[\l] at \p {$\q$};}
			\qpt{qcone}{(-.6,0)}{above}{2}\qpt{qcone}{(.6,0)}{above}{2}
			\node[align=center] at (0,-.68) {antipodal\\pairing};
		\end{tikzpicture}
		\caption{The projective plane $22\times$ with two cone points of order two. Antipodal points of the dashed circle are identified.}
		\label{fig:quotient-pgg}
	\end{figure}
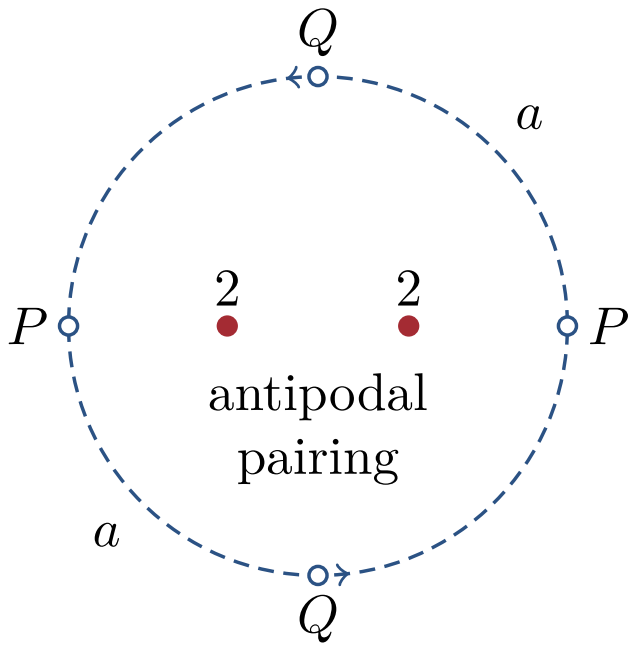
	
	Together, Figures \ref{fig:quotient-smooth}--\ref{fig:quotient-pgg} represent all seventeen rows of Table \ref{quot:tab:types}. The general width formula for every one of these orbifolds is Theorem \ref{thm:affine}.

	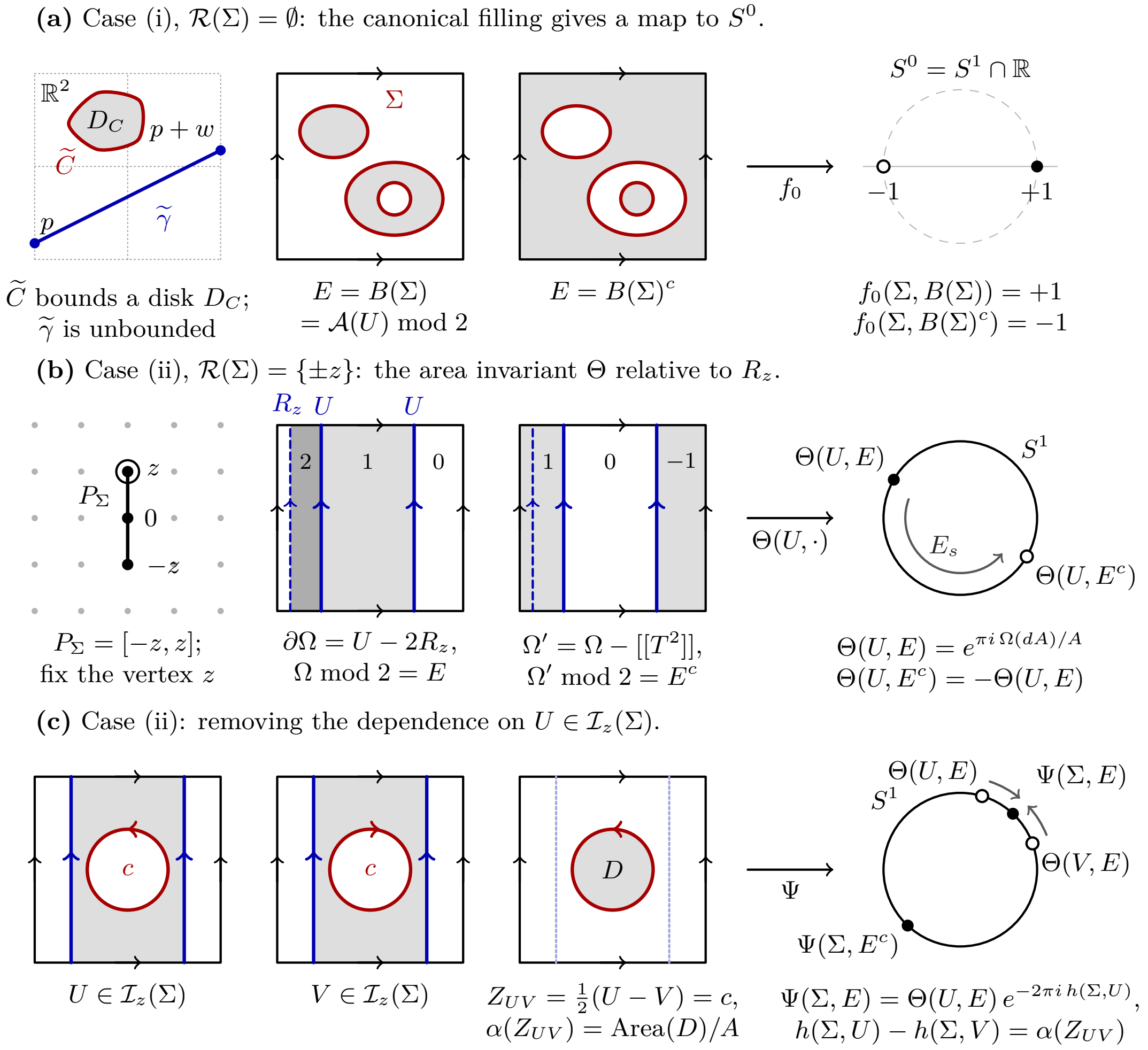
\begin{figure}[tbp]
		\centering
		\colorlet{eqfillone}{black!13}
		\colorlet{eqfilltwo}{black!32}
		\colorlet{eqess}{blue!70!black}
		\colorlet{eqiness}{red!65!black}
		\newcommand{\eqdomain}{%
			\draw[thick] (0,0) rectangle (2.3,2.3);
			\draw[->,thick] (1.0,0) -- (1.32,0);
			\draw[->,thick] (1.0,2.3) -- (1.32,2.3);
			\draw[->,thick] (0,1.0) -- (0,1.32);
			\draw[->,thick] (2.3,1.0) -- (2.3,1.32);}
		\begin{tikzpicture}[font=\footnotesize,line cap=round,line join=round,
			essline/.style={draw=eqess,very thick,postaction={decorate},
				decoration={markings,mark=at position .6 with {\arrow{>}}}},
			refline/.style={draw=eqess,thick,dash pattern=on 2.5pt off 1.6pt,postaction={decorate},
				decoration={markings,mark=at position .6 with {\arrow{>}}}},
			inessline/.style={draw=eqiness,very thick},
			mapsto/.style={->,thick,shorten >=2pt,shorten <=2pt},
			pt/.style={circle,fill=black,inner sep=0pt,minimum size=4.2pt},
			hpt/.style={circle,draw=black,fill=white,thick,inner sep=0pt,minimum size=4.2pt}]
			
			\begin{scope}
				\node[anchor=south west,inner sep=0pt] at (0,2.8)
				{\textbf{(a)} Case (i), $\cR(\Sigma)=\emptyset$: the canonical filling gives a map to $S^0$.};
				\draw[gray!70,densely dotted] (0,0) grid[step=1.15] (2.3,2.3);
				\node[anchor=north west,inner sep=2pt] at (0,2.3) {$\R^2$};
				\filldraw[fill=eqfillone,draw=eqiness,very thick]
				plot[smooth cycle] coordinates {(0.42,1.68) (0.78,2.06) (1.3,1.96) (1.28,1.48) (0.78,1.36)};
				\node at (0.88,1.71) {$D_C$};
				\node[eqiness,anchor=north east,inner sep=1pt] at (0.55,1.42) {$\widetilde C$};
				\draw[eqess,very thick] (0,0.2) -- (2.3,1.35);
				\fill[eqess] (0,0.2) circle[radius=1.8pt] (2.3,1.35) circle[radius=1.8pt];
				\node[anchor=south west,inner sep=1pt] at (0.03,0.27) {$p$};
				\node[anchor=south east,inner sep=1pt] at (2.27,1.42) {$p+w$};
				\node[eqess,anchor=north] at (1.6,0.8) {$\widetilde\gamma$};
				\node[anchor=north,align=center] at (1.15,-0.12)
				{$\widetilde C$ bounds a disk $D_C$;\\ $\widetilde\gamma$ is unbounded};
				\begin{scope}[xshift=3cm]
					\fill[eqfillone,even odd rule]
					(0.7,1.58) ellipse[x radius=.42,y radius=.32]
					(1.45,.75) ellipse[x radius=.6,y radius=.45]
					(1.45,.75) circle[radius=.2];
					\eqdomain
					\draw[inessline] (0.7,1.58) ellipse[x radius=.42,y radius=.32];
					\draw[inessline] (1.45,.75) ellipse[x radius=.6,y radius=.45];
					\draw[inessline] (1.45,.75) circle[radius=.2];
					\node[eqiness] at (1.45,1.98) {$\Sigma$};
					\node[anchor=north,align=center] at (1.15,-0.12) {$E=B(\Sigma)$\\ $\phantom{E}=\cA(U)\bmod 2$};
				\end{scope}
				\begin{scope}[xshift=6cm]
					\fill[eqfillone,even odd rule]
					(0,0) rectangle (2.3,2.3)
					(0.7,1.58) ellipse[x radius=.42,y radius=.32]
					(1.45,.75) ellipse[x radius=.6,y radius=.45]
					(1.45,.75) circle[radius=.2];
					\eqdomain
					\draw[inessline] (0.7,1.58) ellipse[x radius=.42,y radius=.32];
					\draw[inessline] (1.45,.75) ellipse[x radius=.6,y radius=.45];
					\draw[inessline] (1.45,.75) circle[radius=.2];
					\node[anchor=north] at (1.15,-0.12) {$E=B(\Sigma)^c$};
				\end{scope}
				\draw[mapsto] (8.75,1.15) -- node[below] {$f_0$} (9.95,1.15);
				\draw[gray!55,dashed] (11.45,1.15) circle[radius=.95];
				\draw[gray!55] (10.25,1.15) -- (12.65,1.15);
				\node[pt] at (12.4,1.15) {};
				\node[hpt] at (10.5,1.15) {};
				\node[anchor=north] at (12.4,1.08) {$+1$};
				\node[anchor=north] at (10.5,1.08) {$-1$};
				\node[anchor=south] at (11.45,2.12) {$S^0=S^1\cap\R$};
				\node[anchor=north,align=center] at (11.45,-0.12)
				{$f_0(\Sigma,B(\Sigma))=+1$\\ $f_0(\Sigma,B(\Sigma)^c)=-1$};
			\end{scope}
			
			\begin{scope}[yshift=-4.35cm]
				\node[anchor=south west,inner sep=0pt] at (0,2.8)
				{\textbf{(b)} Case (ii), $\cR(\Sigma)=\{\pm z\}$: the area invariant $\Theta$ relative to $R_z$.};
				\foreach \i in {0,...,4}{\foreach \j in {0,...,4}{\fill[gray!60] (0.575*\i,0.575*\j) circle[radius=1.1pt];}}
				\draw[very thick] (1.15,0.575) -- (1.15,1.725);
				\foreach \y in {0.575,1.15,1.725}{\fill (1.15,\y) circle[radius=2pt];}
				\draw[thick] (1.15,1.725) circle[radius=4pt];
				\node[anchor=west] at (1.25,1.75) {$z$};
				\node[anchor=west] at (1.25,0.55) {$-z$};
				\node[anchor=west] at (1.22,1.15) {$0$};
				\node[anchor=east] at (1.08,1.42) {$P_\Sigma$};
				\node[anchor=north,align=center] at (1.15,-0.12) {$P_\Sigma=[-z,z]$;\\ fix the vertex $z$};
				\begin{scope}[xshift=3cm]
					\fill[eqfilltwo] (0.161,0) rectangle (0.544,2.3);
					\fill[eqfillone] (0.544,0) rectangle (1.694,2.3);
					\eqdomain
					\draw[refline] (0.161,0) -- (0.161,2.3);
					\draw[essline] (0.544,0) -- (0.544,2.3);
					\draw[essline] (1.694,0) -- (1.694,2.3);
					\node[eqess,anchor=south] at (0.13,2.3) {$R_z$};
					\node[eqess,anchor=south] at (0.58,2.3) {$U$};
					\node[eqess,anchor=south] at (1.694,2.3) {$U$};
					\foreach \x/\m in {0.3525/2,1.119/1,1.997/0}{\node at (\x,1.85) {\scriptsize$\m$};}
					\node[anchor=north,align=center] at (1.15,-0.12) {$\partial\Omega=U-2R_z$,\\ $\Omega\bmod 2=E$};
				\end{scope}
				\begin{scope}[xshift=6cm]
					\fill[eqfillone] (0,0) rectangle (0.544,2.3);
					\fill[eqfillone] (1.694,0) rectangle (2.3,2.3);
					\eqdomain
					\draw[refline] (0.161,0) -- (0.161,2.3);
					\draw[essline] (0.544,0) -- (0.544,2.3);
					\draw[essline] (1.694,0) -- (1.694,2.3);
					\foreach \x/\m in {0.3525/1,1.119/0,1.997/-1}{\node at (\x,1.85) {\scriptsize$\m$};}
					\node[anchor=north,align=center] at (1.15,-0.12) {$\Omega'=\Omega-[[T^2]]$,\\ $\Omega'\bmod 2=E^c$};
				\end{scope}
				\draw[mapsto] (8.75,1.15) -- node[below] {$\Theta(U,\cdot)$} (9.95,1.15);
				\draw[thick] (11.45,1.15) circle[radius=.95];
				\node[anchor=south west,inner sep=1pt] at ([shift={(45:1.0)}]11.45,1.15) {$S^1$};
				\draw[->,gray!70!black,thick] ([shift={(160:.68)}]11.45,1.15) arc[start angle=160,end angle=320,radius=.68];
				\node[gray!30!black] at ([shift={(240:.4)}]11.45,1.15) {$E_s$};
				\node[pt] (thE) at ([shift={(150:.95)}]11.45,1.15) {};
				\node[hpt] (thEc) at ([shift={(330:.95)}]11.45,1.15) {};
				\node[anchor=south east,inner sep=1pt] at (thE.north west) {$\Theta(U,E)$};
				\node[anchor=north west,inner sep=1pt] at (thEc.south east) {$\Theta(U,E^c)$};
				\node[anchor=north,align=center] at (11.45,-0.12)
				{$\Theta(U,E)=e^{\pi i\,\Omega(dA)/A}$\\ $\Theta(U,E^c)=-\Theta(U,E)$};
			\end{scope}
			
			\begin{scope}[yshift=-8.7cm]
				\node[anchor=south west,inner sep=0pt] at (0,2.8)
				{\textbf{(c)} Case (ii): removing the dependence on $U\in\cI_z(\Sigma)$.};
				\foreach \sh/\name/\tip in {0/U/>,3/V/<}{
					\begin{scope}[xshift=\sh cm]
						\fill[eqfillone,even odd rule] (0.45,0) rectangle (1.85,2.3) (1.15,1.15) circle[radius=.5];
						\eqdomain
						\draw[essline] (0.45,0) -- (0.45,2.3);
						\draw[essline] (1.85,0) -- (1.85,2.3);
						\draw[inessline,postaction={decorate},
						decoration={markings,mark=at position .25 with {\arrow{\tip}}}]
						(1.15,1.15) circle[radius=.5];
						\node[eqiness] at (1.15,1.15) {$c$};
						\node[anchor=north] at (1.15,-0.12) {$\name\in\cI_z(\Sigma)$};
				\end{scope}}
				\begin{scope}[xshift=6cm]
					\eqdomain
					\draw[eqess!35,thick,densely dotted] (0.45,0) -- (0.45,2.3);
					\draw[eqess!35,thick,densely dotted] (1.85,0) -- (1.85,2.3);
					\filldraw[fill=eqfillone,draw=eqiness,very thick,postaction={decorate},
					decoration={markings,mark=at position .25 with {\arrow{>}}}]
					(1.15,1.15) circle[radius=.5];
					\node at (1.15,1.15) {$D$};
					\node[anchor=north,align=center] at (1.15,-0.12)
					{$Z_{UV}=\frac12(U-V)=c$,\\ $\alpha(Z_{UV})=\area(D)/A$};
				\end{scope}
				\draw[mapsto] (8.75,1.15) -- node[below] {$\Psi$} (9.95,1.15);
				\draw[thick] (11.45,1.15) circle[radius=.95];
				\node[anchor=south east,inner sep=1pt] at ([shift={(135:1.0)}]11.45,1.15) {$S^1$};
				\draw[->,gray!70!black,thick] ([shift={(71.5:1.15)}]11.45,1.15) arc[start angle=71.5,end angle=50,radius=1.15];
				\draw[->,gray!70!black,thick] ([shift={(22:1.15)}]11.45,1.15) arc[start angle=22,end angle=43.5,radius=1.15];
				\node[hpt] (thU) at ([shift={(73.5:.95)}]11.45,1.15) {};
				\node[hpt] (thV) at ([shift={(20:.95)}]11.45,1.15) {};
				\node[pt] (psi) at ([shift={(46.75:.95)}]11.45,1.15) {};
				\node[pt] (mpsi) at ([shift={(226.75:.95)}]11.45,1.15) {};
				\node[anchor=south east,inner sep=1pt] at (thU.north) {$\Theta(U,E)$};
				\node[anchor=north west,inner sep=1pt] at (thV.south east) {$\Theta(V,E)$};
				\node[anchor=south west,inner sep=1pt] at ([shift={(47:1.3)}]11.45,1.15) {$\Psi(\Sigma,E)$};
				\node[anchor=north east,inner sep=1pt] at (mpsi.south west) {$\Psi(\Sigma,E^c)$};
				\node[anchor=north,align=center] at (11.45,-0.12)
				{$\Psi(\Sigma,E)=\Theta(U,E)\,e^{-2\pi i\,h(\Sigma,U)}$,\\ $h(\Sigma,U)-h(\Sigma,V)=\alpha(Z_{UV})$};
			\end{scope}
		\end{tikzpicture}
		\caption{The local equivariant maps constructed in Proposition  \ref{prop_equiv_maps}. The (smallest) squares are fundamental domains of $T^2$, red curves are inessential, blue curves are essential, and gray regions are fillings.
			(a) The lift of an inessential loop $C$ to $\R^2$ bounds a disk $D_C$. The lift $\widetilde\gamma$ of an essential closed geodesic is unbounded: the segment from $p$ to $p+w$, $w\in\Lambda$, is one period of $\widetilde\gamma$ and projects to the closed geodesic. The projected disks give the canonical filling $B(\Sigma)$, here a disk and an annulus.
			(b) Two parallel closed geodesics, and we consider two fillings $\Omega$ and $\Omega'$. The two fillings come from currents differing by $[[T^2]]$, so their area invariants are antipodal. Translating $\Sigma$ by half a period (Example \ref{ex:Theta}) exchanges $E$ and $E^c$ and rotates $\Theta$ by $\pi$, so no equivariant map to $S^0$ exists over this family.
			(c) The representatives $U,V$ differ only on $\scD(U,V)=c$, so $\Theta(U,E)/\Theta(V,E)=e^{2\pi i\alpha(Z_{UV})}$. The function $h$ of Lemma  \ref{lem.cocycle} removes this ambiguity; for this single cycle one may take $h(\Sigma,U)=-h(\Sigma,V)=\frac12\alpha(Z_{UV})$, so $\Psi(\Sigma,E)$ is the midpoint of the arc. Angles on the circles are drawn to scale.}
		\label{fig:local-equivariant-maps}
	\end{figure}
	\clearpage

	\bibliographystyle{alpha}
	\bibliography{GMT}
	
\end{document}